\pdfoutput=1
\documentclass[11pt,a4paper]{article}
\usepackage{amsmath,amsthm,amssymb,times}
\usepackage{graphicx}
\usepackage{authblk}
\usepackage{caption, subcaption}
\usepackage[semicolon,round]{natbib}
\usepackage[left=1.5in,top=1.0in,bottom=1.0in,right=1.5in]{geometry}
\numberwithin{equation}{section}
\numberwithin{table}{section}

\newcommand{\dA}{\,\mathrm{d} A}
\newcommand{\intA}{\int_A\!}
\renewcommand{\u}{\mathbf{u}}
\newcommand{\w}{\mathbf{w}}
\newcommand{\F}{\mathbf{F}}
\newcommand{\Q}{\mathbf{Q}}
\renewcommand{\P}[1]{{\mathrm{P}}_{#1}}
\newcommand{\DG}[1]{{\mathrm{P}}_{#1}^\mathrm{DG}}
\newcommand{\pp}[2]{\frac{\partial #1}{\partial #2}}
\newcommand{\dd}[2]{\frac{\delta #1}{\delta #2}}
\newcommand{\total}[2]{\frac{\mathrm{d} #1}{\mathrm{d} #2}}
\newcommand{\scurl}{\nabla^\perp} % 'scalar curl' - continuous + discrete
\newcommand{\vcur}{\nabla^\perp\cdot} % use when nothing follows the character
\newcommand{\vcurl}{\nabla^\perp\!\cdot} % else use this
\newcommand{\wgrad}{\widetilde{\nabla}} % 'weak grad' - discrete only
\newcommand{\wcur}{\widetilde{\nabla}^\perp\cdot} % use when nothing follows the character
\newcommand{\wcurl}{\widetilde{\nabla}^\perp\!\cdot} % else use this

\begin{document}

\author[1,2,*]{Andrew~T.~T.~McRae}
\author[1]{Colin~J.~Cotter}
\affil[1]{\small{\emph{Department of Aeronautics, Imperial College
  London, London, SW7 2AZ, UK}}}
\affil[2]{\small{\emph{The Grantham Institute for Climate Change,
  Imperial College London, London, SW7 2AZ, UK}}}
\affil[*]{\small{Correspondence to:
  \texttt{a.mcrae12@imperial.ac.uk}}}
\title{Energy- and enstrophy-conserving schemes for the shallow-water
  equations, based on mimetic finite elements}
\date{}
\maketitle
\begin{abstract}
This paper presents a family of spatial discretisations of the
nonlinear rotating shallow-water equations that conserve both energy
and potential enstrophy. These are based on two-dimensional mixed
finite element methods and hence, unlike some finite difference
methods, do not require an orthogonal grid. Numerical verification
of the aforementioned properties is also provided.
\end{abstract}

\section{Introduction} The quest for scalable, massively parallel
numerical weather prediction models has led to great interest in
extensions of C-grid staggering to more general mesh structures, such
as icosahedral and cubed meshes. There is also increasing interest in
atmosphere and ocean models that allow arbitrary mesh refinement, in
order to facilitate seamless regional modelling within a global model.
C-grid staggering was proposed as a way of preventing spurious
numerical wave propagation that appears on other grid
staggerings~\citep{sadourny1975dynamics, arakawa1977computational};
these spurious waves interfere with geostrophic adjustment processes
in the numerical solution and rapidly degrade predictive skill. It was
known from the beginning that the C-grid staggering admits
natural finite-difference differential operators (div, grad, curl)
that satisfy discrete versions of vector calculus identities
(div--curl = 0; curl--grad = 0). These
identities allow a separation of the irrotational and solenoidal
components of velocity, which play quite different roles in the low
Rossby number regime.

It was also recognised, from experience with incompressible
quasigeostrophic models~\citep{arakawa1966computational}, that
conservation of energy and potential enstrophy are important for
obtaining nonlinear stability of the model without excessive numerical
diffusion. An energy-conserving formulation was provided
in~\citet{sadourny1975dynamics}, and a formulation that conserves both
energy and enstrophy was given in~\citet{arakawa1981potential}. In the
regime of quasigeostrophic turbulence, the shallow-water equations
exhibit a cascade of energy to large scales. On the other hand,
enstrophy cascades to small scales, and so it makes sense to attempt
to dissipate enstrophy at small scales. The Anticipated Potential
Vorticity Method (APVM) was introduced as a closure to represent the
cascade to scales below the grid
width~\citep{sadourny1985parameterization}; for an appropriate choice
of parameters, the APVM is closely related to Lax--Wendroff advection
schemes. The APVM was incorporated into an energy-conserving,
enstrophy-dissipating shallow-water model
in~\citet{arakawa1990energy}; it remains useful to start with an
enstrophy-conserving model and to then introduce an
enstrophy-dissipating term, since one then has complete control over
the enstrophy dynamics in the model. \citet{arakawa1990energy} also
demonstrated how to handle massless layers in this framework, which
become the basis of many isopycnal ocean
models~\citep[for example]{hallberg1996buoyancy}.

On the sphere, the development of C-grid staggerings for grids other
than the usual latitude-longitude grid was guided by the extension of
the C-grid div, grad and curl operators to arbitrary grids by the
mimetic finite difference community, and by the connection with finite
volume methods~\citep{hyman1997natural}. One route towards energy- and
enstrophy-conserving schemes was proposed using Nambu
brackets~\citep{salmon2005general, salmon2007general,
sommer2009conservative, gassmann2008towards}. In general, a key
challenge was the design of reconstruction methods for the Coriolis
term that allowed for steady linear geostrophic modes on the
$f$-plane, without which nonlinear solutions near to geostrophic
balance would spuriously couple with fast gravity waves. Such a
reconstruction was apparent for triangular grids by making use of the
Raviart-Thomas reconstruction~\citep{bonaventura2005analysis}, but
unfortunately the triangular scheme suffers from spurious branches of
inertia-gravity waves that render it
problematic~\citep{danilov2010utility, gassmann2011inspection}.

A suitable reconstruction on hexagonal grids was then provided
in~\citet{thuburn2008numerical} and extended to arbitrary orthogonal
polygonal grids in~\citet{thuburn2009numerical}, and
energy-conserving, enstrophy-dissipating schemes for the nonlinear
shallow-water equations on arbitrary orthogonal grids were introduced
in~\citet{ringler2010unified}. As discussed
in~\citet{staniforth2012horizontal}, the global degree-of-freedom
ratio between velocity and pressure is altered by increasing or
decreasing the number of cell edges. This may lead to spurious mode
branches -- spurious inertia-gravity wave branches are present for
triangles and spurious Rossby mode branches, for hexagons -- so
quadrilaterals are preferred in order to minimise the possibility of
spurious modes. This suggests the cube mesh for modelling on the
sphere. Unfortunately, the orthogonality requirement in the
construction of~\citet{thuburn2009numerical} leads to meshes that
cluster resolution around the cube vertices, which leads to
non-uniform parallel communication requirements. This
led~\citet{thuburn2012framework} to extend the framework
of~\citet{thuburn2009numerical} to non-orthogonal grids. It has since
been discovered that the scheme of~\citet{thuburn2009numerical} on the
dual icosahedral grid and the scheme of~\citet{thuburn2012framework}
on the cube grid both have inconsistent discretisations of the
Coriolis term (Thuburn, personal communication), meaning that grid
refinement does not improve the accuracy of this term. This, together
with the additional flexibility to alter degree-of-freedom ratios and
to increase the order of accuracy, has motivated the investigation of
mixed finite element methods.

Mixed finite element methods are the analogue of staggered grids since
they use different finite element spaces for velocity and pressure.
Many different combinations of finite element spaces have been
examined in the ocean modelling literature~\citep{le2005dispersion,
le2007analysis, rostand2008raviart, le2008analysis, le2009impact,
danilov2008modeling, comblen2010practical, cotter2011numerical,
le2012spurious}. \citet{cotter2012mixed} concentrated on combinations
of spaces that have discrete versions of the div--curl and curl--grad
identities, just like the C-grid. In the numerical analysis
literature, this is referred to as ``finite element exterior
calculus''~\citep{arnold2006finite}. These combinations were shown to
provide all the properties of the C-grid staggering, including steady
linear geostrophic modes on the $f$-plane, and hence merited further
investigation~\citep{cotter2012mixed}. \citet{staniforth2012analysis}
examined wave propagation for one particular combination, namely the
2nd order Raviart--Thomas ($\mathrm{RT}_1$) space for velocity and
the bilinear discontinuous ($\mathrm{Q}_1^\mathrm{DG}$) space for
pressure, and observed a $2\Delta x$ mode with zero group velocity;
this mode can be corrected by partially lumping the velocity mass
matrix.

In this paper we provide a formulation that uses mixed finite elements
of the type proposed in~\citet{cotter2012mixed}. The formulation
closely follows the steps of~\citet{ringler2010unified}: the
prognostic variables are velocity and layer depth, but there is a
diagnostic potential vorticity that satisfies a discrete conservation
law. Using this potential vorticity in the vector-invariant form of
the equations (as used in the classical C-grid development) naturally
leads to an energy- and enstrophy-conserving form of the equations
without further modification. The conservation properties arise from
the mimetic properties combined with the integral formulation. We
introduce a finite element version of the APVM that dissipates
enstrophy at the gridscale. This formulation is illustrated through
numerical experiments that demonstrate the energy and enstrophy
properties, and demonstrate that the numerical scheme is convergent
and stable. The analytic shallow-water equations and a selection of
derived results are given in section \ref{sec:ana}. We give our
proposed spatial discretisation in section \ref{sec:fed}. Numerical
validation is presented in section \ref{sec:num}, and further areas of
research are discussed in the conclusion. We close by demonstrating
that the conservation properties arise from an almost-Poisson
structure of the spatially discretised equations; this is in Appendix
\ref{sec:poisson}.

\section{Analytic Formulation}
\label{sec:ana}
In this section, we review conservation properties of the rotating
shallow-water equations, since their proofs will be extended to the
finite element discretisations in section \ref{sec:fed}.

The nonlinear shallow-water equations in a rotating frame of reference
are commonly written as \begin{align}
\label{eq:sw-pre1}
\pp{\u}{t} + (\u\cdot\nabla)\u + f\u^\perp &= -g\nabla h\ ,\\
\label{eq:sw-pre2}
\pp{h}{t} + \nabla\cdot(h\u) &= 0\ ,
\end{align}
where ${\u(x,y,t)}$ is the velocity, ${h(x,y,t)}$ is the layer depth,
${f(x,y)}$ is the Coriolis parameter, and $g$ is the gravitational
acceleration. We introduce the $^\perp$ notation for brevity: for
a two-dimensional vector $\w$ in the ${x\text{-}y}$ plane,
${\w^\perp = \hat{\mathbf{z}}\times\w}$, a 90$^\circ$
counterclockwise rotation. If $\w$ is a vector field, this is done pointwise.
We will also use the notation ${\scurl}$
and ${\vcur}$: writing $\nabla$ in components as
${(\partial_x, \partial_y)}$, we have ${\scurl=(-\partial_y,
\partial_x)}$. If $\gamma$ is a scalar field,
\begin{equation}
\scurl\gamma=\left(-\pp{\gamma}{y},\pp{\gamma}{x}\right)\ .
\end{equation}
For a vector field $\w$, with ${\w\equiv(u, v)}$ in components, 
\begin{equation}
\vcurl\mathbf{w}= \pp{v}{x}-\pp{u}{y}\ ,
\end{equation}
a two-dimensional form of $\nabla\times$.

When rewritten in terms of the relative vorticity ${\zeta = \vcurl\u
\equiv \hat{\mathbf{z}}\cdot\nabla\times\u}$, (\ref{eq:sw-pre1}) and
(\ref{eq:sw-pre2}) become
\begin{align}
\label{eq:sw-post1}
\pp{\u}{t} + (\zeta+f)\u^\perp + \nabla\left(gh + \frac{1}{2}|\u|^2\right) &= 0\ ,\\
\label{eq:sw-post2}
\pp{h}{t} + \nabla\cdot(h\u) &= 0\ .
\end{align}
This is the so-called `vector-invariant' form of the equations, which
is the starting point for energy- or enstrophy-conserving formulations
using the C-grid staggering; we shall also use this form here.

We can derive a continuity equation for the absolute vorticity ${\zeta
+ f}$. Defining a potential vorticity ${q = \frac{\zeta + f}{h}}$, we
rewrite (\ref{eq:sw-post1}):
\begin{align}
\label{eq:sw-post3}
\pp{\u}{t} + qh\u^\perp + \nabla\left(gh + \frac{1}{2}|\u|^2\right) &= 0\ .
\intertext{We now apply the ${\vcur}$ operator to (\ref{eq:sw-post3}), giving}
\pp{}{t}(\vcurl\u) + \vcurl(qh\u^\perp)&=0\ , \\
\implies\pp{\zeta}{t} + \nabla\cdot(qh\u)&=0\ .
\end{align}
Assuming ${\pp{f}{t} = 0}$, we then have
\begin{equation}
\label{eq:sw-post4}
\pp{}{t}(qh) + \nabla\cdot(qh\u)=0\ ,
\end{equation}
which is the equation for $q$ written in local conservation form. From
this, we can derive an advection equation for the potential vorticity
$q$. Recall the continuity equation (\ref{eq:sw-post2}). Multiplying this
by $q$, and comparing with (\ref{eq:sw-post4}), we obtain
\begin{equation}
\label{eq:sw-post5}
h\left[\pp{q}{t} + (\u\cdot\nabla)q\right]=0\ ,
\end{equation}
implying that $q$ remains constant in a Lagrangian frame moving with
fluid particles. In particular, if $q$ is initially uniform, $q$ will
remain uniform (and constant) for all time.

In a boundary-free domain, several quantities are conserved.
Integrating (\ref{eq:sw-post4}) over the whole domain gives
conservation of the total absolute vorticity ${\int_A\! qh\dA}$. Less
trivially, the total enstrophy ${\intA q^2h\dA}$ and the total energy
${\intA \left[\frac{1}{2}h|\u|^2 + \frac{1}{2}gh^2\right]\dA}$ are
also constant.

The conservation of enstrophy follows from direct manipulation:
\begin{align}
\total{}{t} \intA q^2h\dA &= \intA \left[2q\pp{}{t}(qh) - q^2\pp{h}{t}\right]\dA \\
&= \intA \left[ 2q\nabla\cdot(-qh\u) - q^2\nabla\cdot(-h\u)\right]\dA \\
&= -\intA \nabla\cdot(q^2h\u)\dA \\
&= 0\ ,\nonumber
\end{align}
where we have used (\ref{eq:sw-post4}) and (\ref{eq:sw-post2}) between
the first and second line. A similar result for higher order moments
of potential vorticity can be obtained by replacing $q^2$ with $q^m$.

Similarly, conservation of energy follows from
\begin{align}
&\total{}{t} \intA \left[\frac{1}{2}h|\u|^2 + \frac{1}{2}gh^2\right]\dA \nonumber \\
&\qquad\qquad= \intA \left[h\pp{}{t}\left(\frac{1}{2}|\u|^2\right)
  + \frac{1}{2}|\u|^2\pp{h}{t} + gh\pp{h}{t}\right]\dA \\
&\qquad\qquad= \int_A\! \left[h\pp{}{t}\left(\frac{1}{2}|\u|^2\right)
  + \left(\frac{1}{2}|\u|^2 + gh\right)\pp{h}{t}\right]\dA \\
&\qquad\qquad= \intA \left[-h\u\cdot\nabla\left(\frac{1}{2}|\u|^2 + gh\right)
  - \left(\frac{1}{2}|\u|^2 + gh\right)\nabla\cdot(h\u)\right]\dA \\
&\qquad\qquad= -\intA \nabla\cdot \left[h\u\left(\frac{1}{2}|\u|^2 + gh\right)\right]\dA \\
&\qquad\qquad= 0\ ,\nonumber
\end{align}
where we have used $\u\cdot$(\ref{eq:sw-post1}) and
(\ref{eq:sw-post2}) between the third and fourth line.

\section{Finite Element Discretisation}
\label{sec:fed}
In this section, we present a family of spatial discretisations, based
on the Finite Element Method, for the nonlinear rotating shallow-water
equations. These discretisations will mimic many properties of the
continuous equations, including the conservation of enstrophy and
energy. The prognostic variables will be the velocity field $\u$ and
the layer depth $h$. Our method explicitly defines a
potential vorticity field $q$ and a volume flux $\F$. However, these
should be interpreted as diagnostic functions of $\u$ and $h$, rather
than independent variables in their own right.

The critical step is the choice of function spaces in which our fields
will reside. In the Finite Element Method, the domain is partitioned
into a large number of non-overlapping subdomains (elements). The
function space specification can be divided into two parts: the
behaviour of a function within each element, and the continuity of a
function at the element boundaries. Almost all function spaces are
piecewise-polynomial (that is, a polynomial when restricted to a
single element). For a scalar function space, the most common
continuity constraints are:
\begin{itemize}
  \item $C^0$ continuous - giving the Continuous Galerkin family
    $\P{n}$, where $n$ is the polynomial degree, and
  \item discontinuous - giving the Discontinuous Galerkin family
    $\DG{n}$, where $n$ is the polynomial degree.
\end{itemize}
Other, less common, conditions include $C^1$ continuity between
elements, and nonconforming ($C^0$ continuity at only the midpoint of
edges). The Continuous and Discontinuous Galerkin families are somewhat
natural function spaces for scalar fields; this can be stated more
precisely in the context of finite element exterior
calculus~\citep{arnold2006finite}. Commonly-used vector function
spaces are often merely tensor products of these two types of scalar
function spaces. However, a careless choice of function space can lead
to genuinely incorrect results, such as spurious solutions arising in
eigenvalue problems~\citep{arnold2010finite}.

We now introduce the function spaces that we will use, and the
relations between them; further details can be found
in~\citet{cotter2012mixed}. We make use of a family of partially
discontinuous vector spaces which are contained in H(div), in other
words they are `div-conforming':
\begin{equation}
\intA \u\cdot\u + (\nabla\cdot\u)(\nabla\cdot\u)\dA < \infty\ .
\end{equation}
Since the functions will be piecewise-polynomial, this condition can
only be violated due to behaviour at element boundaries. The normal
component of the vector field must therefore be continuous across
element boundaries, although the tangential component may be
discontinuous (there is a related space H(curl) in which the opposite
is true). Our velocity field $\u$ and volume flux $\F$ will live in
this space, which we will denote $\mathrm{S}$. Examples include the
Raviart--Thomas family $\mathrm{RT}_n$~\citep{raviart1977mixed}, the
Brezzi--Douglas--Marini family $\mathrm{BDM}_n$~\citep{brezzi1985two},
and the Brezzi--Douglas--Fortin--Marini family
$\mathrm{BDFM}_n$~\citep{brezzi1991mixed}.

For each choice of $\mathrm{S}$, we can define a scalar function space
\begin{equation}
\mathrm{V} = \left\{\nabla\cdot\w:\w\in\mathrm{S}\right\}\ .
\end{equation}
This space is totally discontinuous at element boundaries. The layer
depth $h$ will be in $\mathrm{V}$. Finally, following principles of
Finite Element Exterior Calculus, we define a function space
$\mathrm{E}$ such that
\begin{equation}
\scurl\mathrm{E} \equiv \{\scurl\gamma \colon \gamma\in\mathrm{E}\} \subset \mathrm{S}\ ,
\end{equation}
and
\begin{equation}
\scurl\mathrm{E} = \left\{\ker(\nabla\cdot\colon\mathrm{S}\to\mathrm{V})\right\}\ ;
\end{equation}
$\scurl$ maps bijectively from $\mathrm{E}$ to ${\{\ker(\nabla\cdot)\}\subset\mathrm{S}}$,
modulo constant functions.
This ensures that, for any ${\gamma\in\mathrm{E}}$, ${\nabla\cdot\scurl\gamma\equiv 0}$, the
zero-function in $\mathrm{V}$, and is the analogue of the continuous identity
${\nabla\cdot\nabla\times\equiv 0}$. $\mathrm{E}$ is continuous at element boundaries,
and will contain the potential vorticity field $q$.

We refer to the ${\scurl\colon\mathrm{E}\to\mathrm{S}}$ and ${\nabla\cdot\colon\mathrm{S}\to\mathrm{V}}$
operators as `strong' derivatives, since they act in a pointwise sense
and are identical to the `continuous' $\scurl$ and $\nabla\cdot$ operators.
There are corresponding `weak' operators ${\wcur\colon\mathrm{S}\to\mathrm{E}}$ and
${\wgrad\colon\mathrm{V}\to\mathrm{S}}$ which do not act pointwise, but are instead
defined via integration by parts. Before we elaborate, we take the
opportunity to introduce some notation. We will use angle brackets to
denote the standard $L^2$ inner product:
\begin{equation}
\left\langle f,g \right\rangle = \intA f(x')g(x')\dA\ ,
\qquad\left\langle \u,\mathbf{v} \right\rangle = \intA \u(x')\cdot\mathbf{v}(x')\dA\ .
\end{equation}
Then, in a domain without boundaries, we define $\wcur$ and $\wgrad$ by
\begin{equation}
\left\langle\gamma,\wcurl\u\right\rangle
  = -\left\langle\scurl\gamma,\u\right\rangle\ ,
  \quad\forall\gamma\in\mathrm{E}\ ,
\end{equation}
\begin{equation}
\left\langle \w,\wgrad h\right\rangle
  = -\left\langle\nabla\cdot\w,h\right\rangle\ ,
  \quad\forall\w\in\mathrm{S}\ .
\end{equation}
This is a surprisingly natural definition: let $\Pi_\mathrm{E},
\Pi_\mathrm{S}, \Pi_\mathrm{V}$ be operators that $L^2$-project
arbitrary functions into $\mathrm{E}$, $\mathrm{S}$ and $\mathrm{V}$
respectively, i.e.
\begin{equation}
\left\langle \gamma, \Pi_\mathrm{E}(f)\right\rangle
  = \left\langle \gamma, f\right\rangle\ ,
  \quad\forall \gamma \in \mathrm{E}\ ,
\end{equation}
with $\Pi_\mathrm{S}$ and $\Pi_\mathrm{V}$ defined analogously. Then the following
identities hold:
\begin{equation}
\wcurl(\Pi_\mathrm{S}(\mathbf{v})) \equiv \Pi_\mathrm{E}(\vcurl\mathbf{v})\ ,
\end{equation}
\begin{equation}
\wgrad(\Pi_\mathrm{V}(f)) \equiv \Pi_\mathrm{S}(\nabla f)\ ,
\end{equation}
where $\mathbf{v}$ and $f$ are arbitrary functions; the weak
differential operators commute with $L^2$ projection into the function
spaces. These identities underlie the proof of steady linear
geostrophic modes in the $f$-plane outlined in~\citet{cotter2012mixed}. 

We are now ready to present the discretisation. Recall that $\u$ is
in $\mathrm{S}$ and $h$ is in $\mathrm{V}$. The continuous potential vorticity $q$
satisfied ${qh = \zeta + f}$, where ${\zeta = \vcurl\u}$. However, this
is invalid in our discrete framework: for ${\u\in\mathrm{S}}$, ${\vcurl\u}$ is not
generally defined, since the tangential component of $\u$ is not
continuous across element boundaries. Instead, we must use the weak
operator $\wcur$ discussed previously. Our discrete potential
vorticity ${q \in \mathrm{E}}$ is therefore defined to satisfy, in a
boundary-free domain,
\begin{equation}
\label{eq:fe-vor}
\left\langle \gamma,qh \right\rangle
  = \left\langle -\scurl\gamma, \u \right\rangle
  + \left\langle \gamma,f \right\rangle ,
  \quad\forall \gamma \in \mathrm{E}\ .
\end{equation}
In a domain with boundaries, we would pick up a non-vanishing surface
integral when integrating by parts.

In the continuity equation (\ref{eq:sw-post2}), there was a term
${\nabla\cdot(h\u)}$. Since ${h \in \mathrm{V}}$ is discontinuous, this
expression is problematic. In order to write a discrete continuity
equation, we define a discrete volume flux $\F$ to be the $L^2$
projection of $h\u$ into $\mathrm{S}$, i.e.
\begin{equation}
\label{eq:fe-vol}
\left\langle \w, \F \right\rangle
  = \left\langle \w, h\u \right\rangle\ ,
  \quad\forall \w \in \mathrm{S}\ .
\end{equation}
We can then replace ${\nabla\cdot(h\u)}$ by ${\nabla\cdot\F}$.

Similarly, in the momentum equation (\ref{eq:sw-post3}), there was a
term ${\nabla\left(gh + \frac{1}{2}|\u|^2\right)}$, which is again
incompatible with our discrete framework. We replace $\nabla$ by the
weak gradient $\wgrad$ discussed previously. The discrete forms of our
evolution equations (\ref{eq:sw-post1}) and (\ref{eq:sw-post2}) are
then
\begin{align}
\label{eq:fe-evo1}
\left\langle \w,\pp{\u}{t} \right\rangle
  + \left\langle \w,q\F^\perp\right\rangle
  - \left\langle \nabla\cdot\w, gh+\frac{1}{2}|\u|^2 \right\rangle &= 0\ ,
  \quad\forall \w \in \mathrm{S}\ ,\\
\label{eq:fe-evo2}
\left\langle \phi, \pp{h}{t} \right\rangle
  + \left\langle \phi,\nabla\cdot\F \right\rangle &= 0\ ,
  \quad\forall \phi \in \mathrm{V}\ .
\end{align}
The equations (\ref{eq:fe-vor}) through (\ref{eq:fe-evo2}) form our
scheme. Note that (\ref{eq:fe-evo1}) holds even in a domain with
boundaries, as long as ${\u\cdot\mathbf{n}=0}$. More importantly, 
(\ref{eq:fe-evo2}) implies that the equation
\begin{equation}
\label{eq:fe-mass}
\pp{h}{t} + \nabla\cdot\F = 0
\end{equation}
is satisfied pointwise, as both $\pp{h}{t}$ and $\nabla\cdot\F$ are in $\mathrm{V}$.

In a boundary-free domain, these discrete equations reproduce the
results given in the previous section for the continuous governing
equations. Recalling that ${\nabla\cdot\scurl\gamma\equiv 0}$, we begin by
inserting ${\w = -\scurl\gamma}$ into (\ref{eq:fe-evo1}), for any
${\gamma \in \mathrm{E}}$:
\begin{equation}
\left\langle -\scurl\gamma,\pp{\u}{t} \right\rangle
  + \left\langle -\scurl\gamma,q\F^\perp\right\rangle =0\ ,
  \quad\forall\gamma\in\mathrm{E}\ .
\end{equation}
Assuming that ${\pp{f}{t} = 0}$, we can rewrite the first
  term using $\pp{}{t}$(\ref{eq:fe-vor}):
\begin{align}
\label{eq:fe-pv1}
\left\langle \gamma,\pp{}{t}(qh) \right\rangle
  + \left\langle -\scurl\gamma,q\F^\perp\right\rangle &=0\ ,
  \quad\forall\gamma\in\mathrm{E} \\
\label{eq:fe-pv2}
\implies \left\langle \gamma,\pp{}{t}(qh) \right\rangle
  + \left\langle -\nabla\gamma,q\F\right\rangle &= 0\ ,
  \quad\forall\gamma\in\mathrm{E} \\
\label{eq:fe-pv3}
\implies \left\langle \gamma,\pp{}{t}(qh) \right\rangle
  + \left\langle \gamma,\nabla\cdot(q\F)\right\rangle &= 0\ ,
  \quad\forall\gamma\in\mathrm{E}\ ,
\end{align}
where the integration by parts in the final line is permitted,
i.e.\ it is an exact identity, since $\gamma$ is continuous and
$\F$ is div-conforming. (\ref{eq:fe-pv3}) is a discrete approximation
to the local conservation law for $q$ (\ref{eq:sw-post4}), which was
previously combined with the continuity equation to form an advection
equation or $q$ (\ref{eq:sw-post5}). A similar procedure can be
carried out in the discrete case by expanding out the derivatives:
\begin{equation}
\left\langle \gamma,h\pp{q}{t} + q\pp{h}{t} \right\rangle
  + \left\langle \gamma,q\nabla\cdot\F + (\F\cdot\nabla) q\right\rangle = 0,
  \quad\forall\gamma\in\mathrm{E}\ .
\end{equation}
We now use our observation (\ref{eq:fe-mass}), which stated that the
continuity equation holds pointwise, implying
\begin{equation}
\label{eq:fe-pv1a}
\left\langle \gamma,h\pp{q}{t} \right\rangle
  + \left\langle \gamma, (\F\cdot\nabla) q\right\rangle = 0,
  \quad\forall\gamma\in\mathrm{E}\ .
\end{equation}
This is a discrete analogue of (\ref{eq:sw-post5}), and is enough to
reproduce the result that if $q$ is initially constant, $q$ remains
constant for all time.

To reproduce conservation laws, we will typically make a specific
choice of the `test-function' $\gamma$ (or $\w$, or $\phi$). For
example, taking ${\gamma \equiv 1}$ in (\ref{eq:fe-pv1}) or
(\ref{eq:fe-pv2}) gives conservation of absolute vorticity in a
boundary-free domain.

Conservation of enstrophy follows from choosing ${\gamma=q}$ (which is
permitted since ${q\in E}$):
\begin{align}
\total{}{t} \intA q^2h\dA
  &\equiv \total{}{t} \left\langle q,qh \right\rangle \\
&= 2\left\langle q,\pp{}{t}(qh)\right\rangle
  - \left\langle q^2,\pp{h}{t}\right\rangle\ . \\
\intertext{Using our result from (\ref{eq:fe-mass}), that
  ${\pp{h}{t} + \nabla\cdot\F = 0}$ is satisfied pointwise, and
  taking ${\gamma = q}$ in (\ref{eq:fe-pv2}):}
&= 2\left\langle \nabla q,q\F\right\rangle
  + \left\langle q^2,\nabla\cdot\F \right\rangle \\
&= \intA \nabla\cdot(q^2\F)\dA \\
&=0\ .\nonumber
\end{align}

Conservation of energy is again obtained by direct computation:
\begin{align}
&\total{}{t} \intA \left[\frac{1}{2}h|\u|^2 + \frac{1}{2}gh^2\right]\dA 
  \equiv \total{}{t}\left(\frac{1}{2}\left\langle h\u,\u\right\rangle
  + \frac{1}{2}\left\langle gh, h\right\rangle\right) \\
&\qquad\qquad= \left\langle h\u,\pp{\u}{t} \right\rangle
  + \left\langle \pp{h}{t},\frac{1}{2}|\u|^2\right\rangle
  + \left\langle \pp{h}{t},gh\right\rangle\ . \\
\intertext{Using (\ref{eq:fe-vol}) with ${\w=\pp{\u}{t}}$ (permitted since
  ${\pp{\u}{t}\in S}$), we obtain}
&\qquad\qquad= \left\langle \F,\pp{\u}{t}\right\rangle
  + \left\langle \pp{h}{t},gh + \frac{1}{2}|\u|^2\right\rangle. \\
\intertext{Then, using (\ref{eq:fe-evo1}) with ${\w=\F}$ (permitted since
  ${\F\in S}$), and (\ref{eq:fe-mass}), we obtain}
&\qquad\qquad= \left\langle \F,-q\F^\perp \right\rangle
  + \left\langle \nabla\cdot\F, gh+\frac{1}{2}|\u|^2 \right\rangle
  + \left\langle -\nabla\cdot\F, gh + \frac{1}{2}|\u|^2 \right\rangle \\
&\qquad\qquad= 0\ ,\nonumber
\end{align}
as required. An explanation of how these properties arise from a
discrete almost-Poisson structure is provided in Appendix
\ref{sec:poisson}.

Equations (\ref{eq:fe-vor}) through (\ref{eq:fe-evo2}) imply a set of
ODEs in the basis coefficients for $\u$ and $h$, which can then be
integrated using any chosen time integration scheme. For explicit
schemes, they will still require the solution of matrix-vector systems
in order to obtain $\pp{\u}{t}$ and $\pp{h}{t}$; the matrices are,
however, very well-conditioned (the condition number being independent
of mesh resolution~\citep{wathen1987realistic}) and, in the case of
$h$, block diagonal.

There is a problem, though: \eqref{eq:fe-pv1a} is the usual Galerkin
finite element discretisation of the advection equation, which, just
like the centred finite difference discretisation, is known to be
unstable when used with explicit time integration
methods~\citep{gresho1998incompressible}. This means that the $L^2$
norm of $q$ will grow without bound, implying that $\u$ will become
increasingly rough.  Additionally, for low Rossby number solutions of
the shallow-water equations near to geostrophic balance, enstrophy is
known to cascade to small scales. This means that an
enstrophy-conserving scheme will lead to a pile up of enstrophy at
small scales, and it is necessary to dissipate enstrophy at such
scales in order to obtain physical solutions. This is an identical
situation to the energy- and enstrophy-conserving scheme
of~\citet{arakawa1981potential}, and indeed to any
enstrophy-conserving scheme. To obtain a stable scheme, we must make
modifications so that equation \eqref{eq:fe-pv2} takes the form
\begin{equation}
\left\langle \gamma, \pp{}{t}(qh)\right\rangle + \langle -\nabla\gamma, q\F
+\Q^*
\rangle = 0,
\end{equation}
where $\Q^*$ is an additional numerical flux that leads to stability --
necessary for convergence of numerical solutions.  This
changes the evolution equation (\ref{eq:fe-evo1}) to the following:
\begin{equation}
\label{eq:fe-evo1-dissipative}
\begin{split}
&\left\langle \w,\pp{\u}{t} \right\rangle
  + \left\langle \w,q\F^\perp+(\Q^*)^\perp\right\rangle \\
&\qquad - \left\langle \nabla\cdot\w, gh+\frac{1}{2}|\u|^2 \right\rangle = 0\ ,
  \quad\forall \w \in \mathrm{S}.
\end{split}
\end{equation}
If, in addition, the dissipative flux $\Q^*$ is proportional to $\F$,
the energy is still conserved, since the corresponding term vanishes
in equation \eqref{eq:fe-evo1-dissipative} when $\w=\F$. The term
$\Q^*$ is chosen so that the divergence-free component of $\u$ remains
stable. In the low Rossby number limit near to geostrophic balance,
the irrotational component of $\u$ is extremely weak and it is not
necessary to introduce further stabilisation to control that
component. Since $\Q^*$ is introduced to dissipate instabilities
generated by the advection term in the PV equation, it evolves on the
slow timescale and therefore does not create a strong source of
inertia-gravity waves; it instead just modifies the ``slow manifold''
about which the fast waves oscillate.

There are a wide range of higher-order time integration schemes
available for the advection equation using continuous finite element
spaces, many of which can be written in the form of the addition of a
dissipative flux $\Q^*$ to discrete counterparts of
\eqref{eq:fe-evo1-dissipative}, including SUPG
\citep{brooks1982streamline} and Taylor-Galerkin methods
\citep{donea1984taylor}. To ease the exposition in this paper by
avoiding complicated discussion of time-discretisation methods and to
provide a link with the history of the development of C-grid grids,
here, following~\citet{arakawa1990energy}, we will introduce the
Anticipated Potential Vorticity
Method~\citep{sadourny1985parameterization} to stabilise the scheme,
by setting $\Q^*=-\tau(\u\cdot\nabla)q\F$ in the continuous time
equations, where $\tau$ is a timescale. By design, this dissipates
enstrophy at small scales by using an upwinded $q$ value in the
advective term, while the conservation of energy is unchanged. The
other equations remain unchanged.  Since we
are using the APVM purely for stabilisation, rather than as a subgrid
parameterisation, we will simply take ${\tau = \frac{\Delta t}{2}}$.
This means that when we discretise the equations in time,
the resulting numerical scheme will be first-order accurate in time.

\section{Numerical results}
\label{sec:num}
The above equations were integrated using the classical 4th order
Runge-Kutta scheme, making use of tools from the FEniCS project: a
collection of free software for automated and efficient solutions of
differential equations~\citep{logg2012automated}. In particular we
make use of the H(div) elements (in this case, $\mathrm{RT}_0$,
$\mathrm{BDM}_1$, $\mathrm{BDM}_2$, and $\mathrm{BDFM}_1$) whose
implementation in FEniCS is described in \cite{rognes2009efficient}.
The goal of the numerical experiments is to demonstrate: a) that they
produce convergent discretisations of the equations, b) that the
claimed energy and enstrophy conservation properties hold, and c) that
they reproduce convincing vortex dynamics within balanced
solutions. All the integrations were performed in planar geometries.

The analytic results derived in the previous section hold for any
function spaces $\mathrm{E}$, $\mathrm{S}$ and $\mathrm{V}$ satisfying
the stated relationships. In this section we will explicitly use
the four triples
${\left(\P{1}, \mathrm{RT}_0, \P{0}\right)}$, ${\left(\P{2},
\mathrm{BDM}_1, \P{0}\right)}$, ${\left(\P{2}\oplus\mathcal{B}_3,
\mathrm{BDFM}_1, \DG{1}\right)}$ and ${\left(\P{3}, \mathrm{BDM}_2,
\DG{1}\right)}$, which adhere to the criteria.

\begin{figure}
  \begin{subfigure}[b]{0.45\columnwidth}
    \includegraphics[width=\columnwidth]{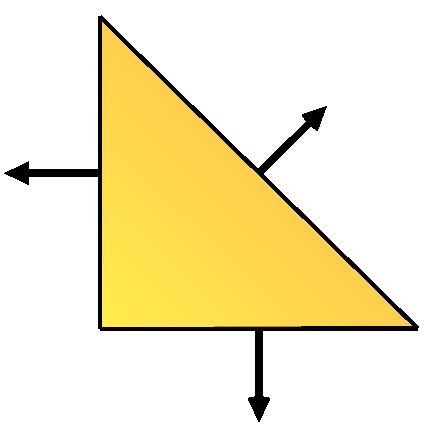}
    \caption{$\mathrm{RT}_0$}
  \end{subfigure}
  \begin{subfigure}[b]{0.45\columnwidth}
    \includegraphics[width=\columnwidth]{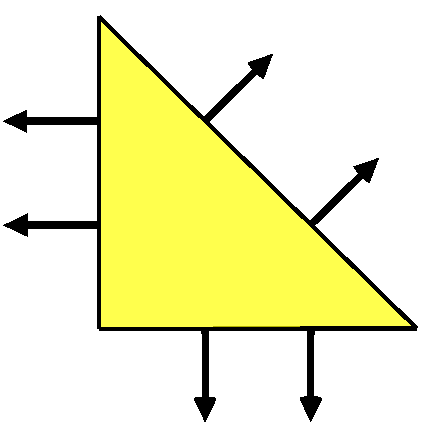}
    \caption{$\mathrm{BDM}_1$}
  \end{subfigure}
  \\
  \begin{subfigure}[b]{0.45\columnwidth}
    \includegraphics[width=\columnwidth]{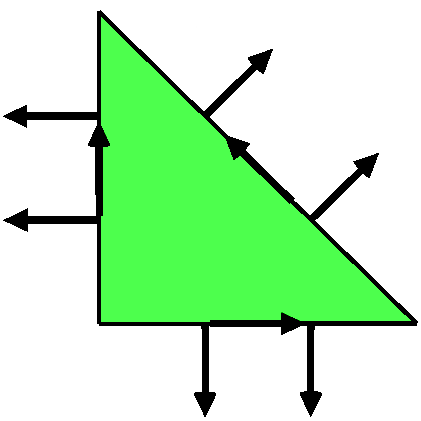}
    \caption{$\mathrm{BDFM}_1$}
  \end{subfigure}
  \begin{subfigure}[b]{0.45\columnwidth}
    \includegraphics[width=\columnwidth]{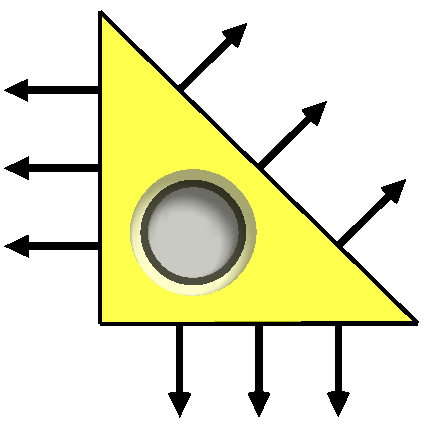}
    \caption{$\mathrm{BDM}_2$}
  \end{subfigure}
\caption{The degrees of freedom for the different velocity function spaces.
$\mathrm{RT}_0$ requires the zeroth moment of the normal component on edges
or, equivalently, point evaluation of the normal component at the midpoint
of each edge. $\mathrm{BDM}_1$ requires zeroth and first moments on edges,
or two point evaluations. $\mathrm{BDFM}_1$ additionally requires the zeroth
moment of tangential velocity on each edge, local to each cell, since
the tangential velocity can be discontinuous between neighbouring cells.
Finally, $\mathrm{BDM}_2$ requires three pointwise evaluations of normal
velocity on each edge, plus three additional interior moments.}
\label{fig:functionspaces}
\end{figure}

The $\P{n}$ and $\DG{n}$ spaces have been introduced already, in the
previous section. $\mathrm{RT}_n, \mathrm{BDM}_n$ and $\mathrm{BDFM}_n$
are the Raviart--Thomas, Brezzi--Douglas--Marini and
Brezzi--Douglas--Fortin--Marini families
respectively~\citep{raviart1977mixed,brezzi1985two,brezzi1991mixed},
and the $n$ suffix indicates a spatial discretisation of order $n+1$.
These somewhat-uncommon vector-valued function spaces are shown in
Figure~\ref{fig:functionspaces}.
${\P{2}\oplus\mathcal{B}_3}$ denotes a continuous, piecewise-quadratic
function enriched with a cubic `bubble' local to each element.

It is known that $\mathrm{RT}$ spaces on triangles have a
surplus of pressure degrees of freedom [DOFs] and consequently have spurious
inertia-gravity modes. $\mathrm{BDM}$ spaces have a deficit of pressure DOFs and
consequently have spurious Rossby modes. $\mathrm{BDFM}_1$ has an exact balance of
velocity and pressure degrees of freedom, which is a necessary
condition for the absence of spurious modes~\citep{cotter2012mixed}, hence
its inclusion in our tests.

Although we will only present results for the four triples mentioned above,
any member of the infinite families ${\left(\P{n}, \mathrm{RT}_{n-1},
\DG{n-1}\right)}$ and ${\left(\P{n+1}, \mathrm{BDM}_n, \DG{n-1}\right)}$
could be used, and three of our four triples are from said families
($\DG{0}$ and $\P{0}$ are synonymous). Also, as discussed in the previous
section, the choice of the velocity space $\mathrm{S}$ determines $\mathrm{V}$
and $\mathrm{E}$. Therefore, from here onwards, we will only refer to the
velocity space used -- $\mathrm{RT}_0$, $\mathrm{BDM}_1$, $\mathrm{BDFM}_1$ or
$\mathrm{BDM}_2$ -- when presenting our results.

\begin{figure}
\minipage{0.45\textwidth}
  \includegraphics[width=\textwidth]{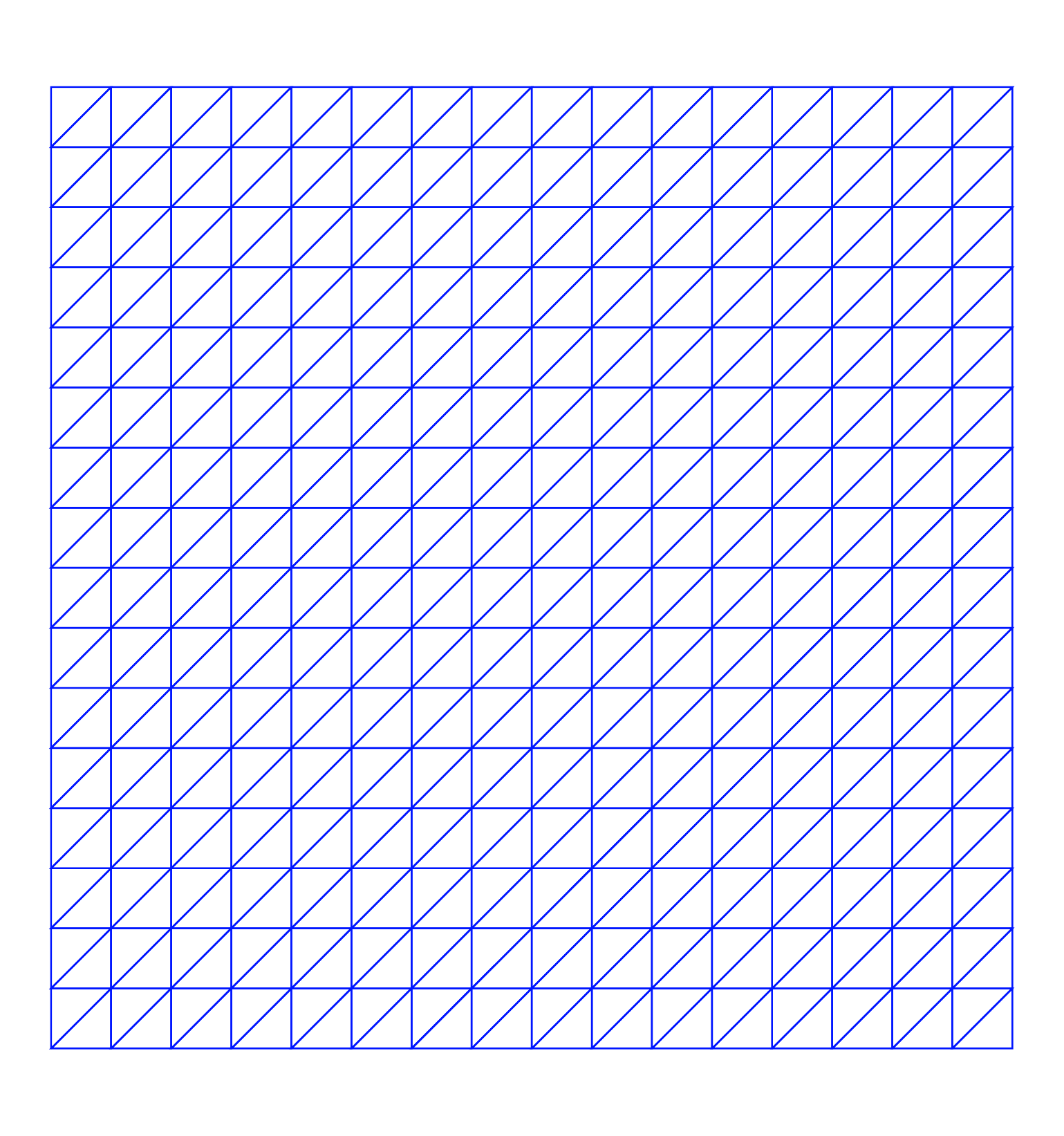}
\endminipage\hfill
\minipage{0.45\textwidth}
  \includegraphics[width=\textwidth]{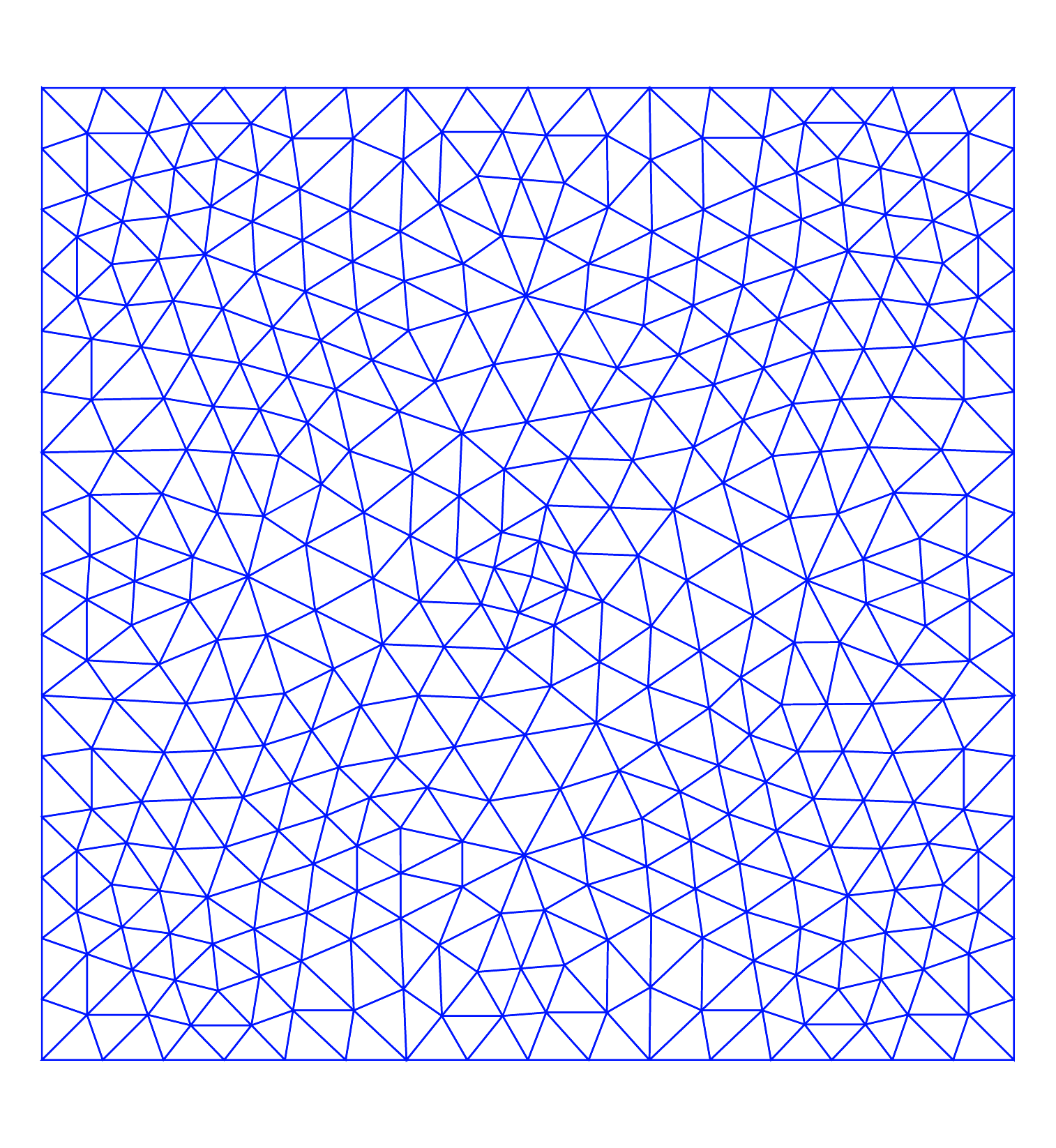}
\endminipage
\caption{Examples of regular and unstructured meshes.}
\label{fig:mesh1}
\end{figure}

To emulate a boundary-free domain, we used $[0, 1]^2$ equipped with
periodic boundary conditions throughout. All lengths are therefore
non-dimensional. We used both regular and unstructured meshes;
examples are given in Figure \ref{fig:mesh1}. The regular meshes are
available in FEniCS by default. The unstructured meshes were generated
in gmsh~\citep{geuzaine2009gmsh} with `target element size'
$\frac{1}{8}$, $\frac{1}{12}$, $\frac{1}{16}$, $\frac{1}{24}$ and
$\frac{1}{32}$.
This gave grids with 160, 416, 736, 1488 and 2744 triangles
respectively. For the unstructured grids, we have plotted errors
against the total number of DOFs. For $\mathrm{RT}_0$,
there are 1.5 global velocity DOFs and 1 height DOF per triangle. For
$\mathrm{BDM}_1$, the corresponding numbers are 3 and 1. For
$\mathrm{BDFM}_1$, 6 and 3; for $\mathrm{BDM}_2$, 7.5 and 3.

We will begin by examining the original, unstabilised scheme, and
verifying that the discrete conservation results indeed hold. We will
then look at the effects of the APVM stabilisation.

\subsection{Balanced state}
\label{sec:balstate}
We performed a convergence test to verify that our implementation is
correct. Here, we restricted ourselves to solutions of the form
${\u = (u(y), 0)}$, ${h = h(y)}$, and $f$ constant. Then the
shallow-water equations reduce to
\begin{equation}
\label{eq:sw-bal}
\frac{\partial{u}}{\partial{t}} = 0,\qquad fu =
  -g\frac{\partial{h}}{\partial{y}},\qquad \pp{h}{t} = 0\ .
\end{equation}
This is a simple example of geostrophic balance, in which the Coriolis
force balances the pressure term exactly, and the advection terms
vanish. 

For our tests we made the particular choice
\begin{equation}
  \begin{split}
  u &= \sin(4\pi y)\ ,\\
  h &= 10 + \frac{1}{4\pi}\frac{f}{g}\cos(4\pi y)\ ,
  \end{split}
\end{equation}
where we have nondimensionalised time accordingly (recall that the domain had
non-dimensional width 1). We will take ${f = 10.0}$ and ${g = 10.0}$, with
the appropriate nondimensionalisations, giving a Rossby number
${Ro \equiv \frac{UL}{f} \simeq 0.1}$ and a Burger number
${B \equiv \frac{gH}{f^2 L^2} \simeq 1}$.
We used RK4 timestepping with ${\Delta t = 0.0005}$ until ${t = 1}$, a
regime in which timestepping error is far smaller than spatial
discretisation error.

\begin{figure}
\centering
\includegraphics[width=0.9\columnwidth]{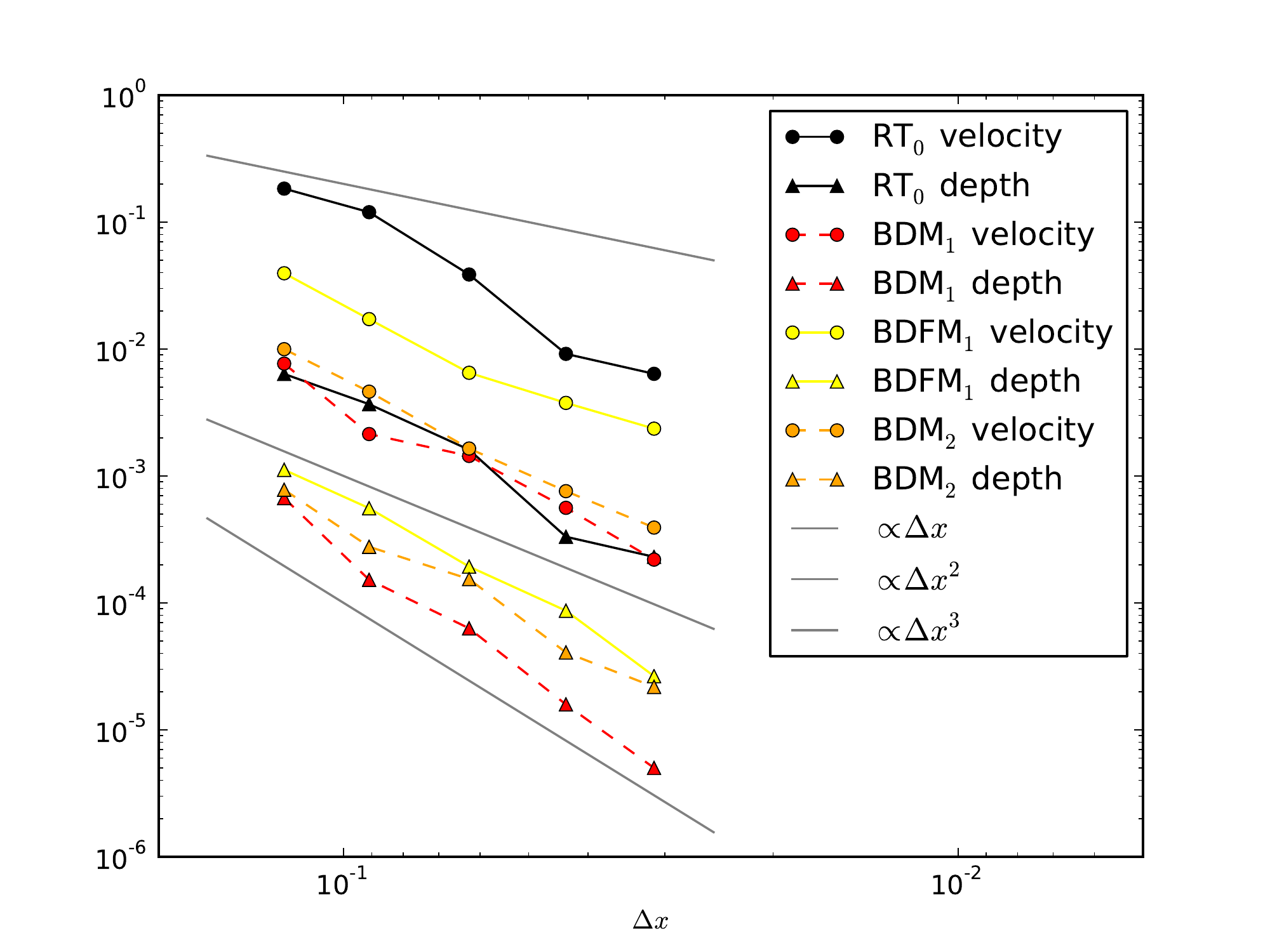}
\caption{$L^2$ norms of relative velocity and height errors when
  simulating the balanced state described in section \ref{sec:balstate}, with
  the unstabilised scheme, on a regular mesh. Error plotted against $\Delta x$.}
\label{fig:balstr-dx}
\end{figure}

\begin{figure}
\centering
\includegraphics[width=0.9\columnwidth]{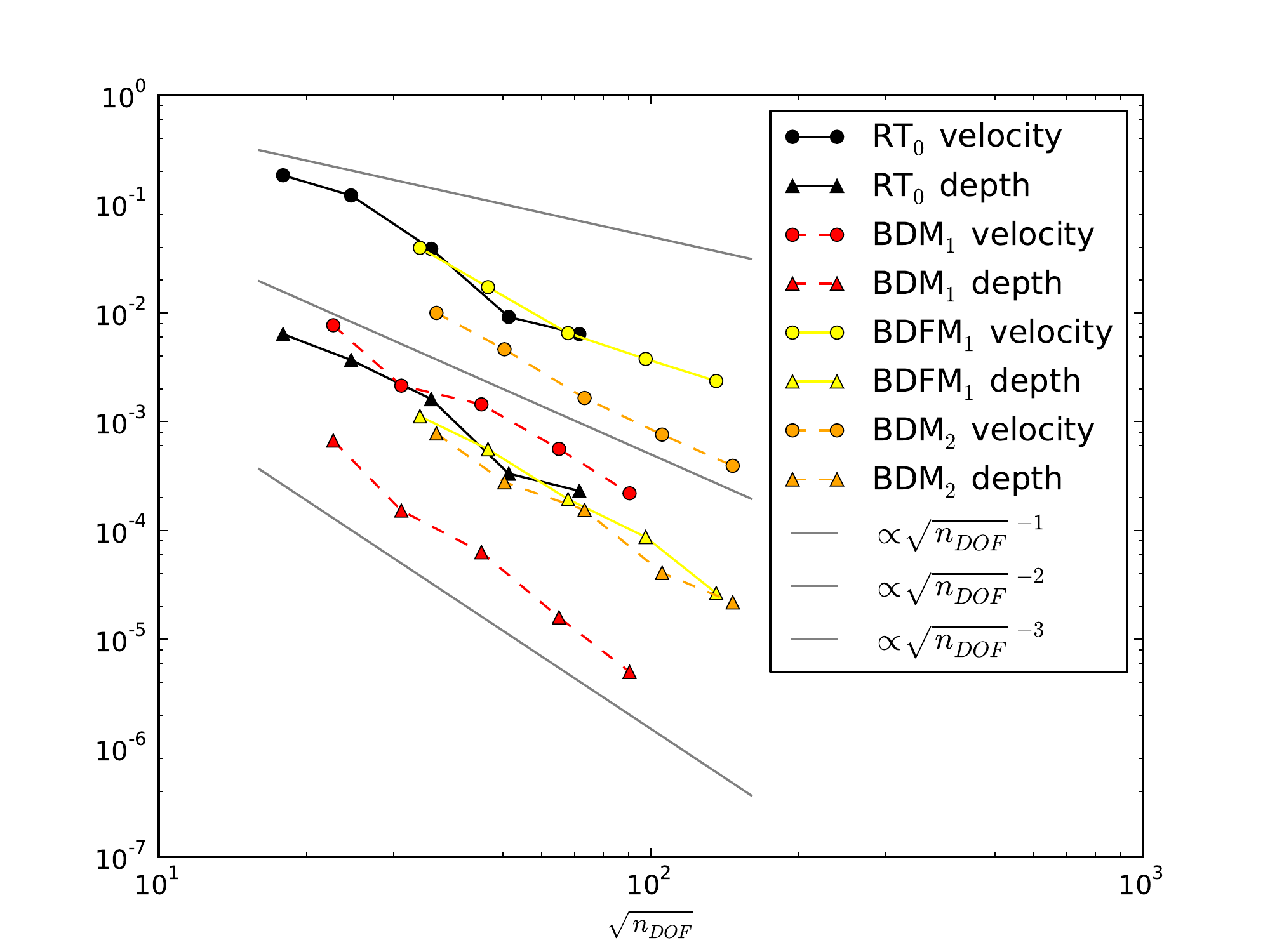}
\caption{$L^2$ norms of relative velocity and height errors when
  simulating the balanced state described in section \ref{sec:balstate},
  with the unstabilised scheme, on a regular mesh. Error plotted against
  the square root of $n_{DOF}$.}
\label{fig:balstr-dof}
\end{figure}

\begin{figure}
\centering
\includegraphics[width=0.9\columnwidth]{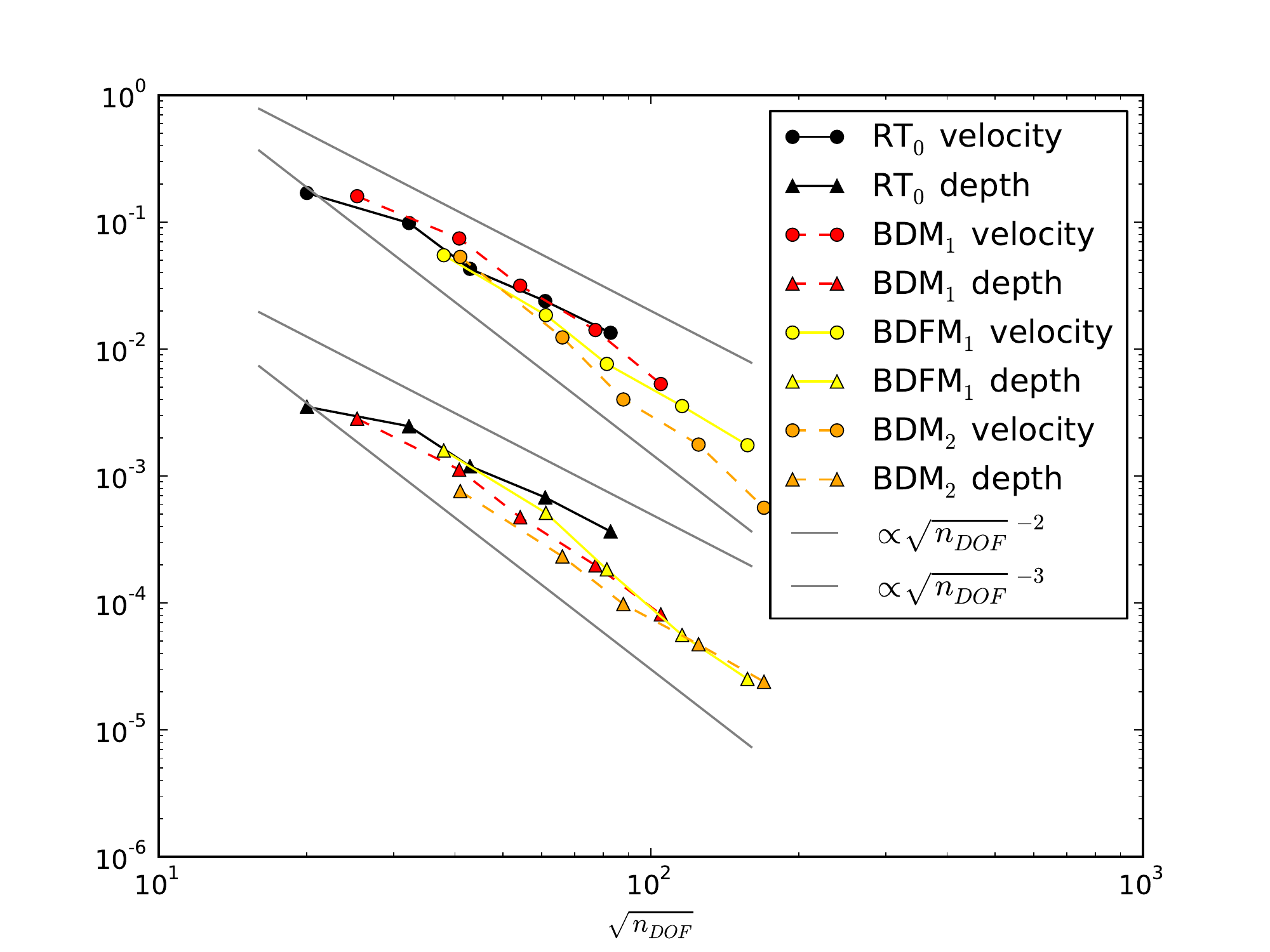}%
\caption{$L^2$ norms of relative velocity and height errors when
  simulating the balanced state described in section \ref{sec:balstate},
  with the unstabilised scheme, on an unstructured mesh.}
\label{fig:balunstr}%
\end{figure}

The $L^2$ norms of ${\u_\textrm{final} - \u_\textrm{initial}}$ and
${h_\textrm{final} - h_\textrm{initial}}$ are shown in figures
\ref{fig:balstr-dx} and \ref{fig:balstr-dof} for a structured
mesh, and figure \ref{fig:balunstr} for an unstructured mesh. We see
at least second-order convergence for all the schemes. This is an
order more than we would naively expect for
$\mathrm{RT}_0$. $\mathrm{BDFM}_1$ and $\mathrm{BDM}_2$ both have
quadratic representations of $q$ which may explain the third order
convergence, which is especially noticeable on the unstructured grid.

\subsection{Energy and Enstrophy conservation}
\label{sec:enens}
\begin{figure}
\centering
\includegraphics[width=0.9\columnwidth]{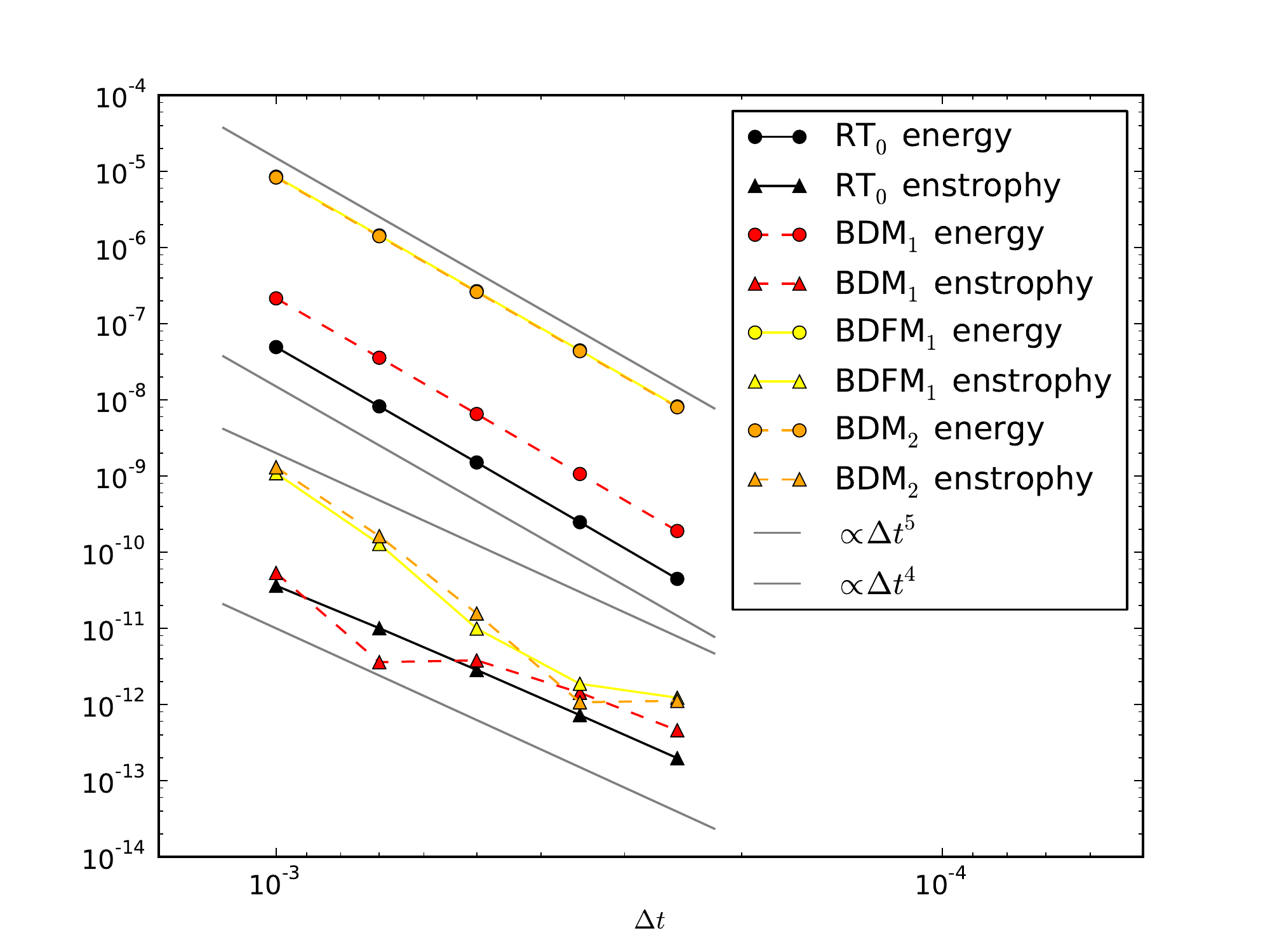}%
\caption{Relative energy and enstrophy errors when the initial condition
  given in section \ref{sec:enens} is simulated, with the unstabilised scheme,
  on a regular mesh with ${\Delta x = \frac{1}{16}}$.}
\label{fig:consstr}%
\end{figure}

\begin{figure}
\centering
\includegraphics[width=0.9\columnwidth]{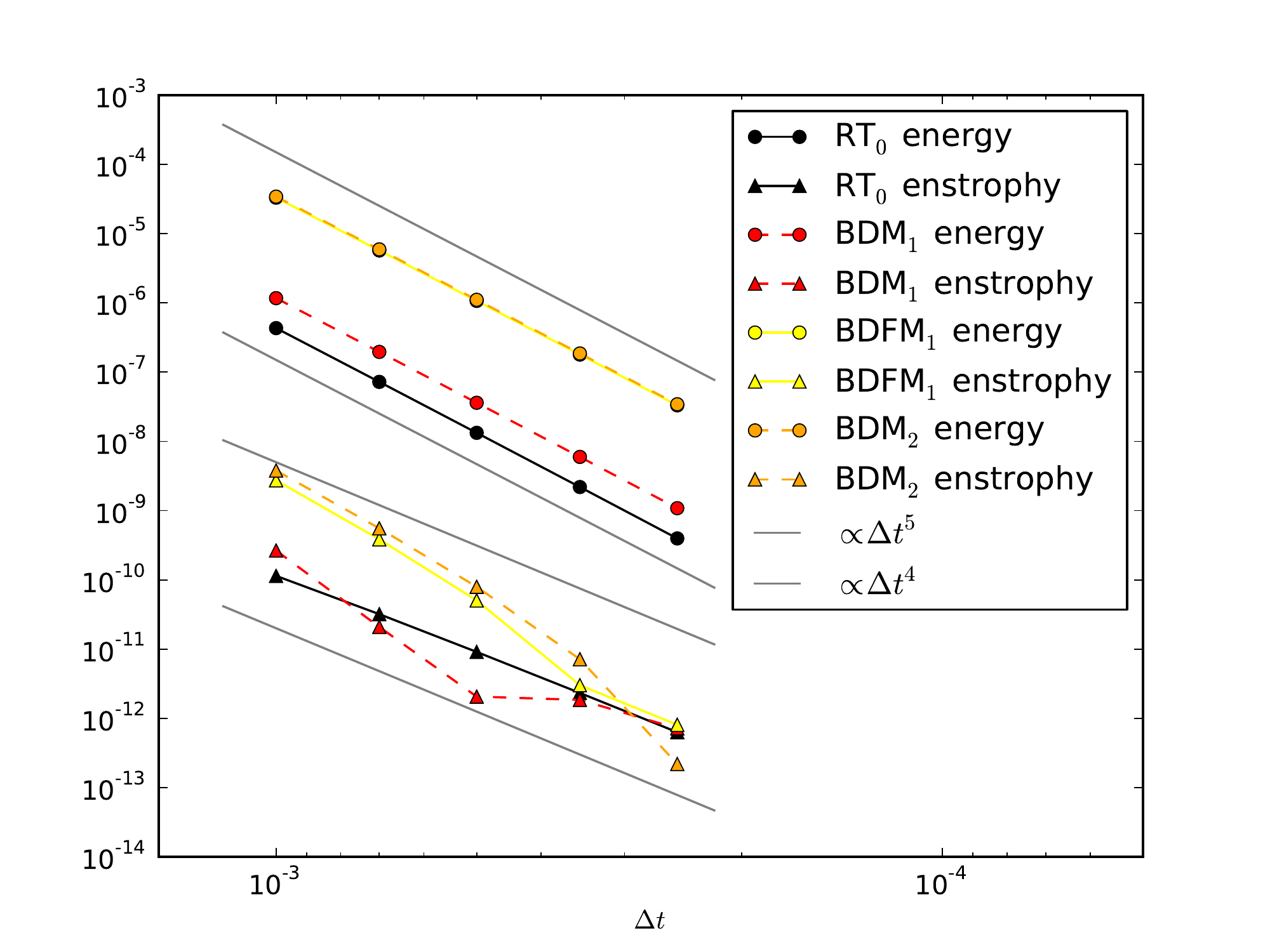}%
\caption{Relative energy and enstrophy errors when the initial condition
  given in section \ref{sec:enens} is simulated, with the unstabilised scheme,
  on an unstructured mesh with 736 triangles.}
\label{fig:consunstr}%
\end{figure}

To demonstrate energy and enstrophy conservation, we took an arbitrary
initial condition and parameters ${f = 5.0}$, ${g = 5.0}$. The system
was simulated with RK4 timestepping for a range of $\Delta t$ until
${t = 1.001}$.  Although the spatial discretisation conserves energy
and enstrophy, the temporal discretisation does not. We expect to see
at most fourth-order errors in the conservation of energy and
enstrophy, with changing $\Delta t$, as the discrete-time numerical
solutions approach the continuous-time, discrete-space solutions. We
used the initial condition
\begin{equation}
  \begin{split}
  \u &= (0,v(x)) = (0,\sin(2\pi x)) \\
h &= h(y) = 1 + \frac{1}{4\pi}\frac{f}{g}\sin(4\pi y)
  \end{split}
\end{equation}
The relative changes in energy and enstrophy between the initial and
final states are shown in figures \ref{fig:consstr} and
\ref{fig:consunstr}. The former is for a regular mesh with ${\Delta x
= \frac{1}{16}}$, the latter for an unstructured mesh with 736
triangles. In both cases, the enstrophy change is fourth-order in
$\Delta t$.  The energy change is fifth-order in $\Delta t$; we believe
that this is due to additional cancellations in the equation for
energy evolution.

\subsection{Stabilised scheme}
\label{sec:stab}
\begin{figure}
\centering
\includegraphics[width=0.9\columnwidth]{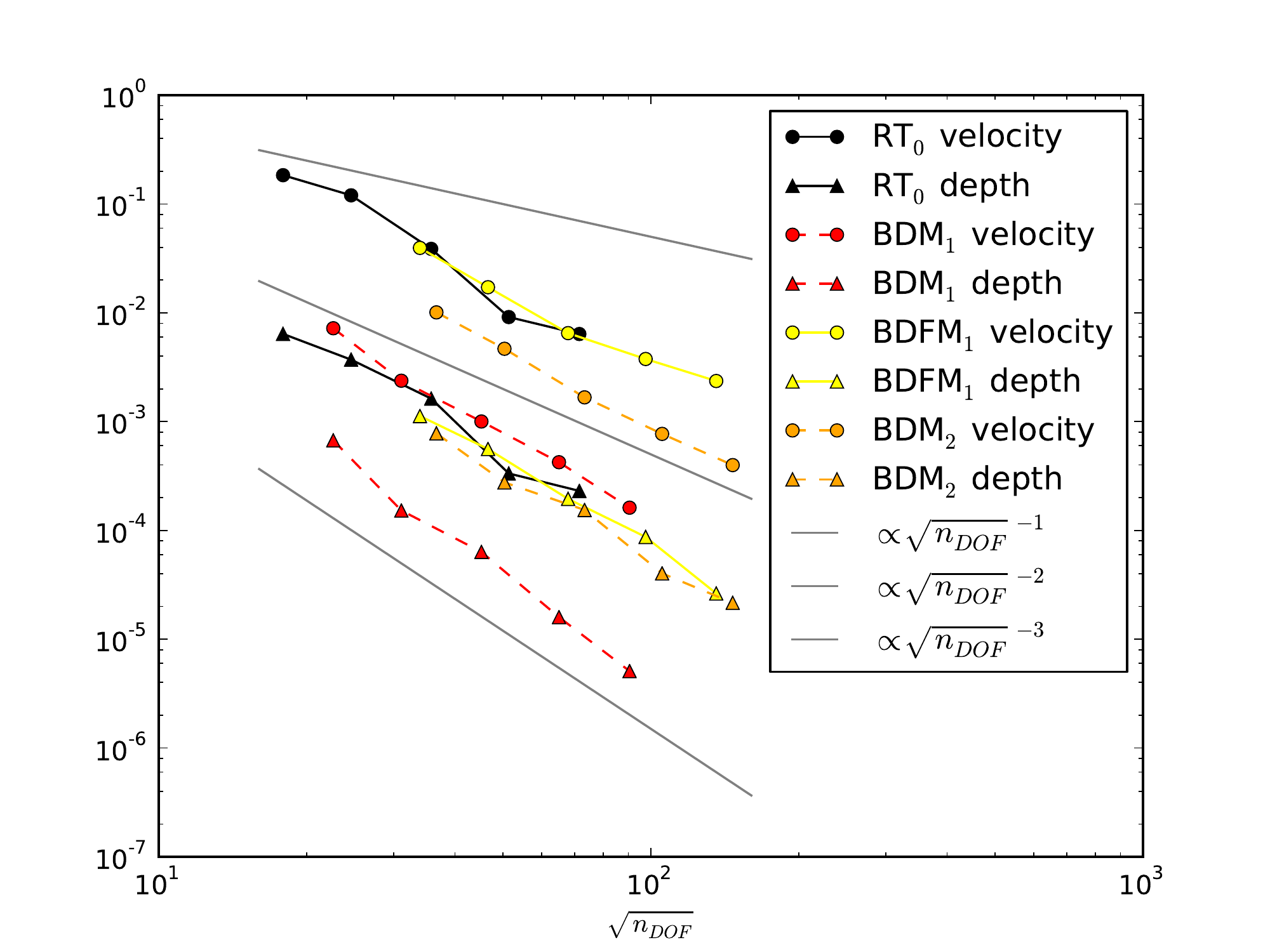}%
\caption{$L^2$ norms of relative velocity and height errors when
  simulating the balanced state described in section \ref{sec:balstate},
  with the stabilised scheme, on a regular mesh.}
\label{fig:APVMbal-str}%
\end{figure}

\begin{figure}
\centering
\includegraphics[width=0.9\columnwidth]{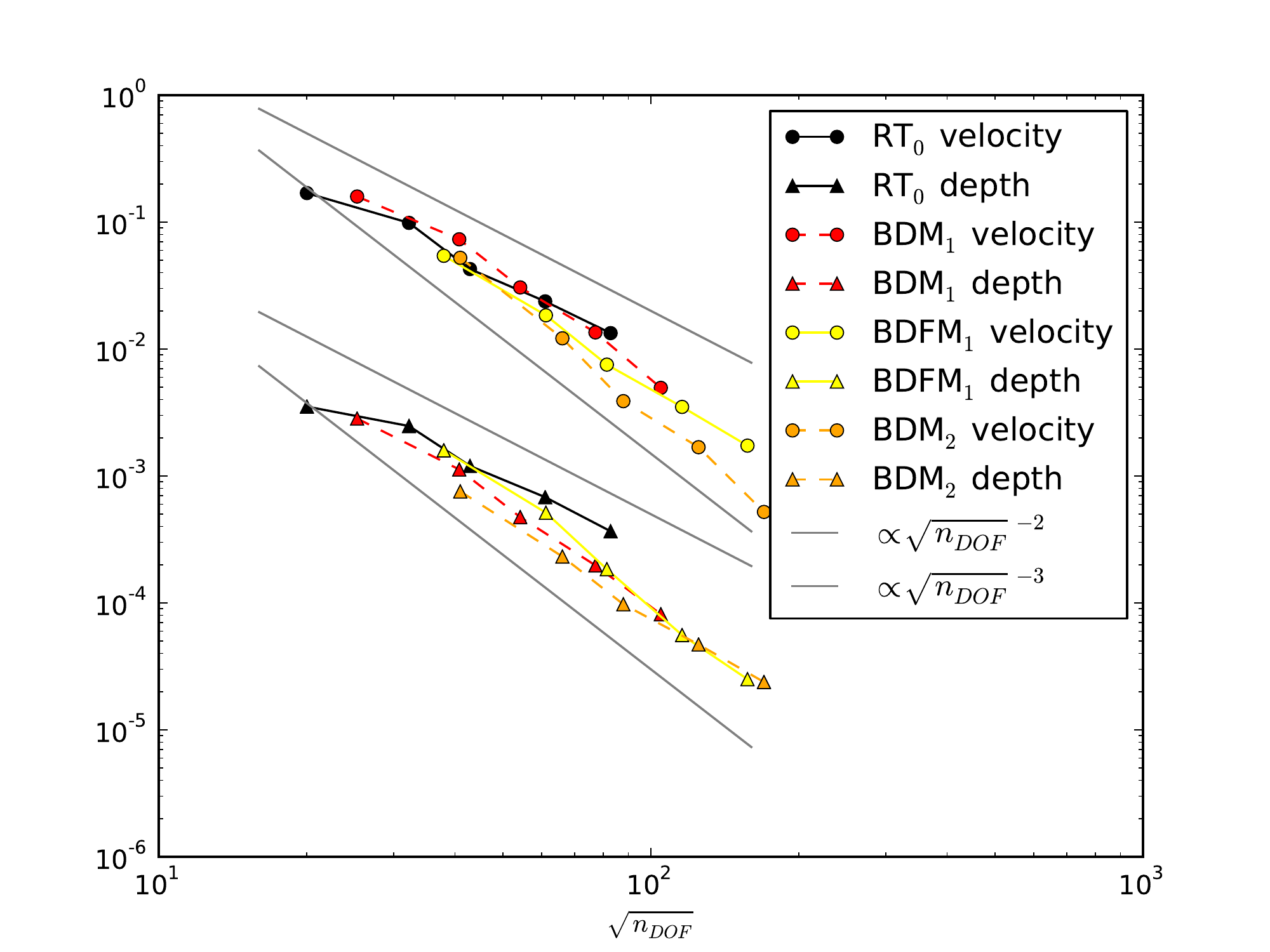}%
\caption{$L^2$ norms of relative velocity and height errors when
  simulating the balanced state described in section \ref{sec:balstate},
  with the stabilised scheme, on an unstructured mesh.}
\label{fig:APVMbal-unstr}%
\end{figure}

We repeated the balanced state convergence test for the scheme
stabilised by the APVM. The $L^2$ norms of ${\u_\textrm{final} -
\u_\textrm{initial}}$ and ${h_\textrm{final} - h_\textrm{initial}}$
are shown in figures \ref{fig:APVMbal-str} and \ref{fig:APVMbal-unstr}
for a regular and unstructured grid, respectively. Note that the
numerical values from the stabilised scheme are almost identical to the
unstabilised scheme, to within a couple of percent.

\begin{figure}
\centering
\includegraphics[width=0.9\columnwidth]{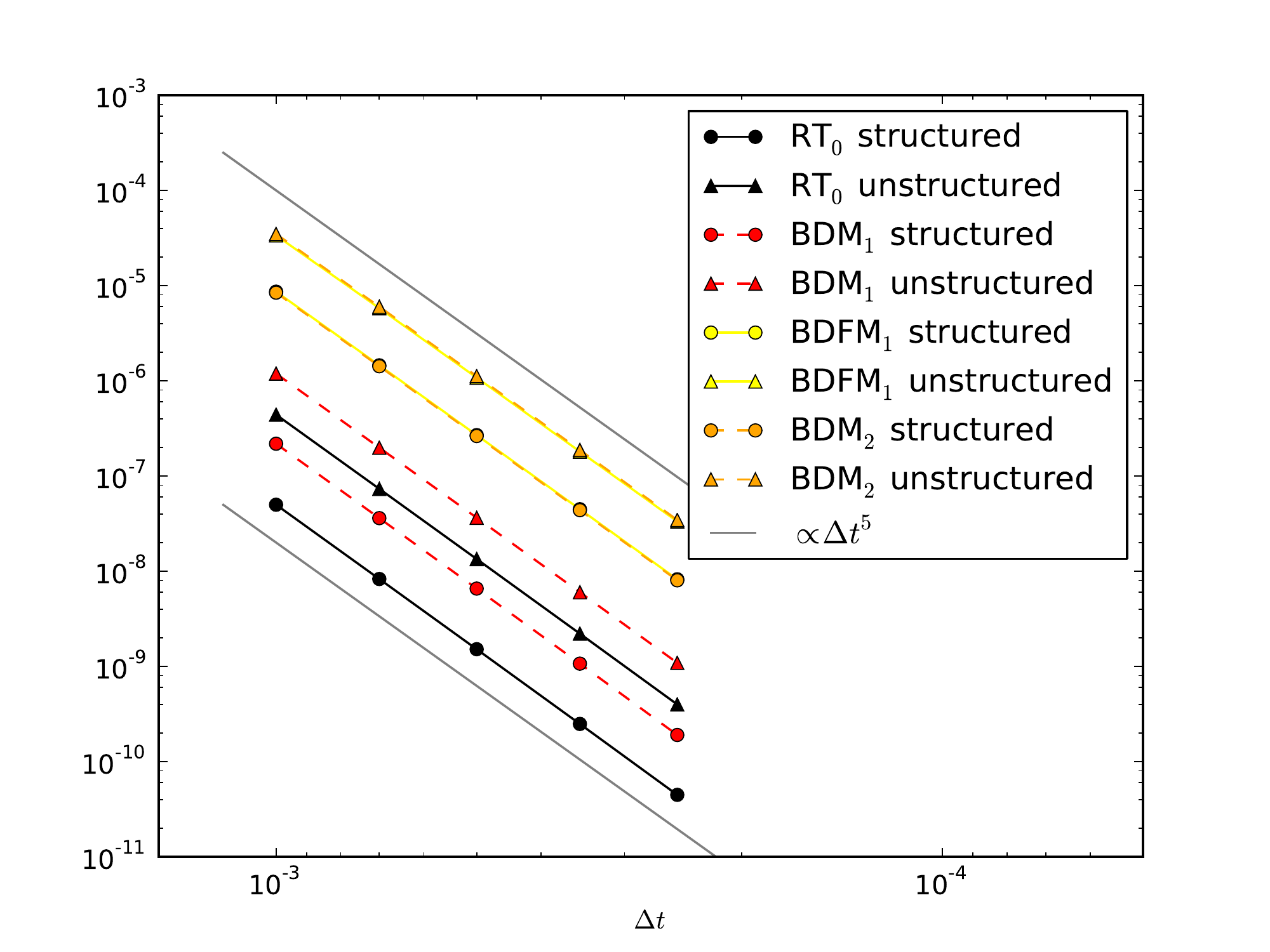}%
\caption{Relative energy error using the Anticipated Potential Vorticity
  Method to stabilise the proposed scheme. As before, it appears to be
  fifth-order in $\Delta t$, consistent with the use of RK4 timestepping.}
\label{fig:APVMeng}%
\end{figure}

\begin{figure}
\centering
\includegraphics[width=0.9\columnwidth]{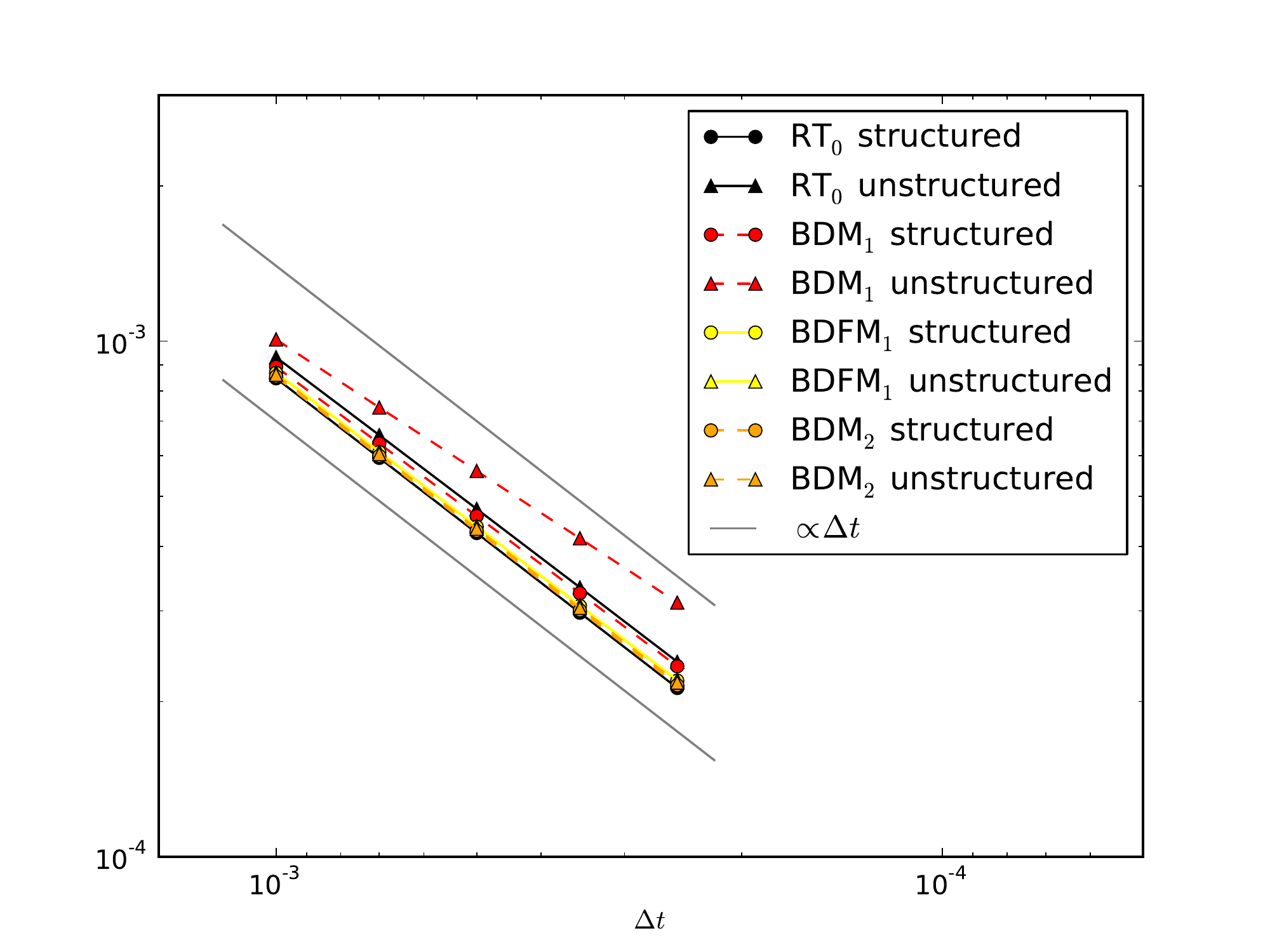}%
\caption{Relative enstrophy change using the Anticipated Potential
  Vorticity Method to stabilise the proposed scheme. As the APVM erodes
  enstrophy, we no longer see fourth-order convergence. First-order
  convergence is seen, since we took ${\tau = \frac{\Delta t}{2}}$}
\label{fig:APVMens}%
\end{figure}

We tested for energy conservation using the same initial conditions as
were used in section \ref{sec:enens}, on the same unstructured grid,
and examined the enstrophy loss. These are shown in figures
\ref{fig:APVMeng} and \ref{fig:APVMens} respectively. As before, the
energy change appears to be at least fourth-order in $\Delta t$ while,
as expected, enstrophy is now dissipated.

\begin{figure}
\centering
\includegraphics[width=0.9\columnwidth]{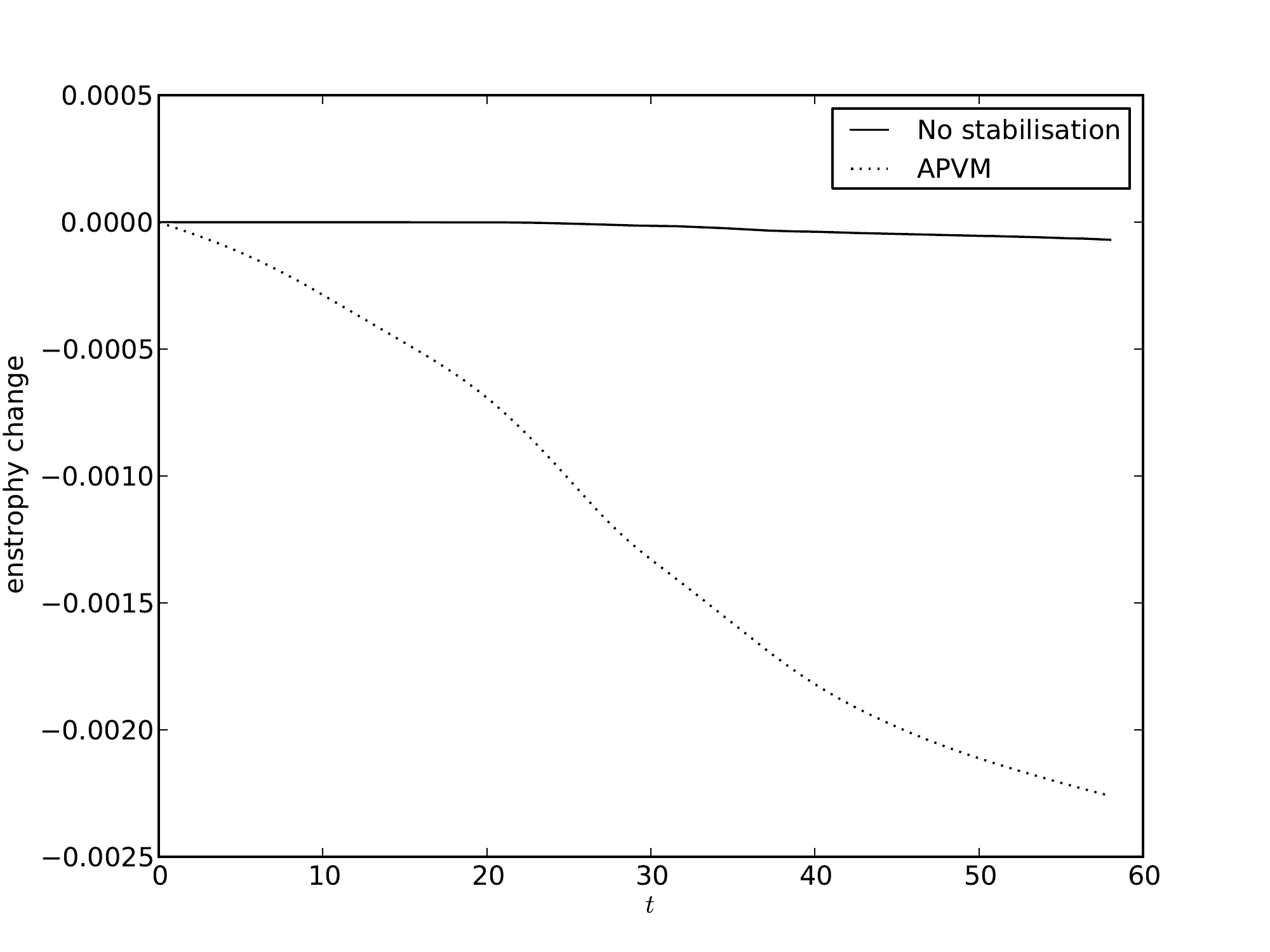}%
\caption{Evolution of total enstrophy in the `merging vortex' problem.
  The stabilised scheme loses a macroscopic amount of enstrophy, while
  the unstabilised scheme only loses enstrophy due to numerical error.
}
\label{fig:ensdecay}%
\end{figure}

\begin{figure}
\centering
  \includegraphics[width=0.7\textwidth]{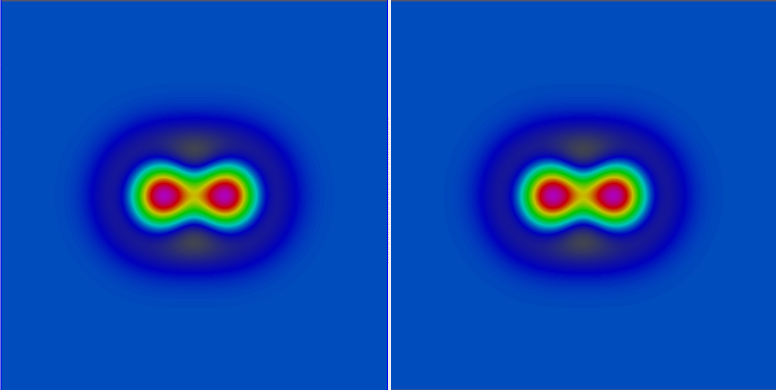}
  \includegraphics[width=0.7\textwidth]{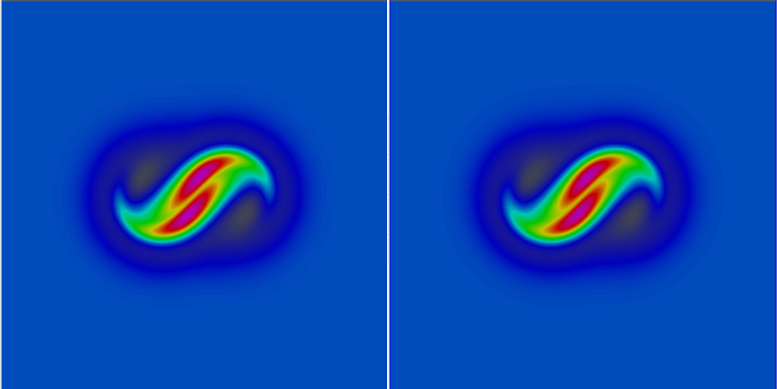}
  \includegraphics[width=0.7\textwidth]{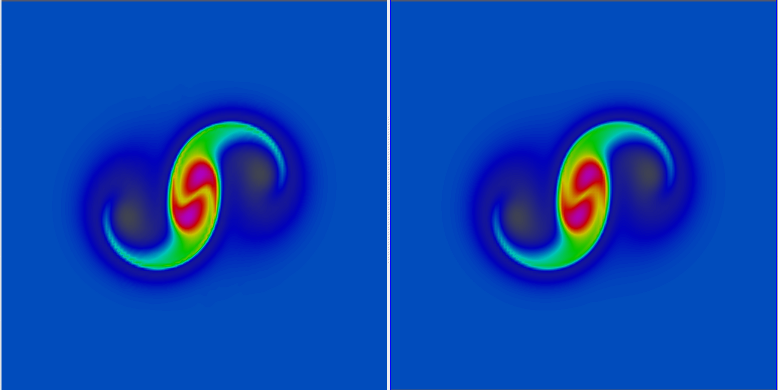}
  \includegraphics[width=0.7\textwidth]{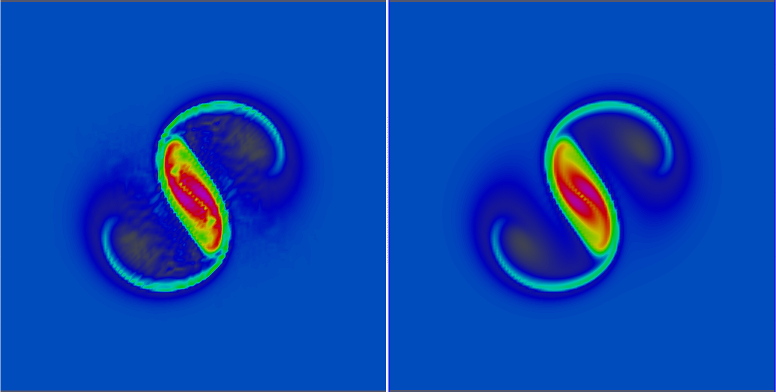}
  \caption{Evolution of merging vortices. The potential vorticity $q$
    is shown, with the stabilised scheme on the right. By the fourth pair
    of images, spurious oscillations are visible when the
    unstabilised scheme is used. The plots above correspond to
    $t$ = 0, 8, 16, 24.}
  \label{fig:vort1}
\end{figure}

\begin{figure}
\centering
  \includegraphics[width=0.7\textwidth]{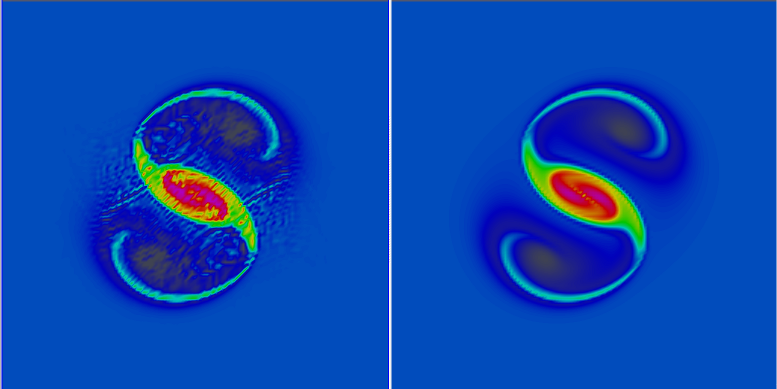}
  \includegraphics[width=0.7\textwidth]{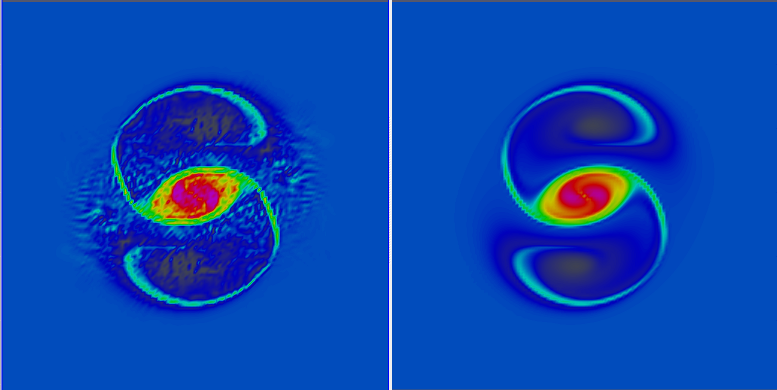}
  \includegraphics[width=0.7\textwidth]{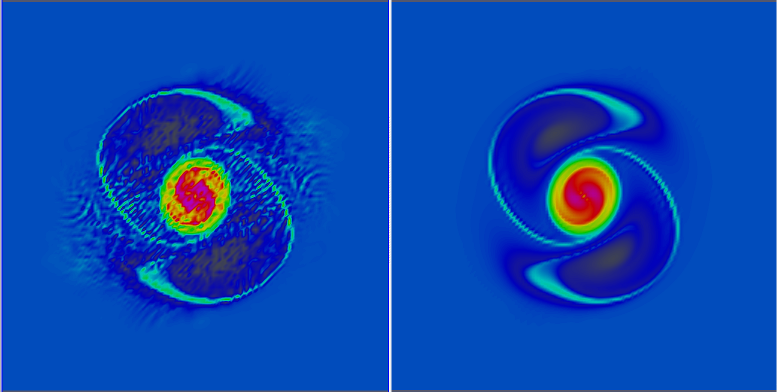}
  \includegraphics[width=0.7\textwidth]{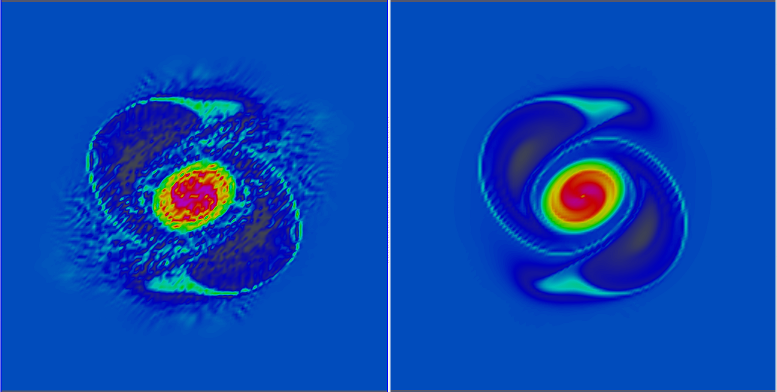}
  \caption{Evolution of merging vortices. The potential vorticity $q$
    is shown, with the stabilised scheme on the right. Spurious
    oscillations are clearly visible when the unstabilised scheme is
    used. The plots above correspond to $t$ = 32, 40, 48, 56.}
  \label{fig:vort2}
\end{figure}

Finally, in figures \ref{fig:vort1} and \ref{fig:vort2}, we show the
evolution of a `merging vortex' problem, in a quasi-geostrophic
parameter regime, in order to visually compare the stabilised and
unstabilised schemes. The initial condition for the velocity field is
derived from a streamfunction: a superposition of two
radially-symmetric Gaussians with different centrepoints. The initial
condition for the depth field is chosen to satisfy linear geostrophic
balance. The $\mathrm{BDM}_1$ function space was used for these
examples. Enstrophy evolution is shown in figure
\ref{fig:ensdecay}. This example demonstrates the ability of the APVM
to dissipate enstrophy on an unstructured mesh in this framework
whilst preserving energy (up to timestepping error). The $L_2$ norm of
the linear geostrophic imbalance $f\u^\perp+g\nabla h$ was calculated
at each timestep, and the differences between with and without APVM
were orders of magnitude smaller than the variation in the imbalance
in either case, which in itself was very small, demonstrating that
APVM does not generate fast inertia-gravity waves.

\clearpage

\clearpage

\clearpage

\section{Conclusion}
In this paper, we introduced a discretisation of the nonlinear
shallow-water equations that extends the energy- and
enstrophy-conserving formulation of~\citet{arakawa1981potential}, and
the energy-conserving, enstrophy-dissipating formulation
of~\citet{arakawa1990energy}, to the mixed finite element approach
advocated in~\citet{cotter2012mixed}. The extension is obtained by
replacing the discrete differential operators defined on the C-grid by
div and curl operators that map between different finite element
spaces. Given these operators, the steps are then identical to the
C-grid approach: a discrete volume flux is obtained, a potential
vorticity is diagnosed and the discrete volume flux is used to create
a discrete potential vorticity flux. This flux is then used in the
vector-invariant form of the equation for $\u$. The energy- and
enstrophy-conservation arises from a discrete Poisson bracket
structure, to be discussed in the Appendix. The convergence and
energy/enstrophy properties of the scheme were demonstrated using
numerical examples.

In ongoing work, we are developing semi-implicit versions of this
discretisation approach, as well as extending it to curved elements
for meshing the sphere, with the aim of prototyping horizontal
discretisations for the UK GungHo Dynamical Core project. We are also
exploring the replacement of \eqref{eq:fe-vol} with an upwind
discontinuous Galerkin scheme (which would dissipate potential energy
at the gridscale) to avoid solution of a global mass matrix, and the
use of explicit Taylor-Galerkin schemes to extend the time accuracy of
the implied PV equation whilst maintaining stability. We are also
investigating the extension of the finite element framework to
three-dimensional flows.

\appendix
  
\section{Almost-Poisson structure of the spatial discretisation}
\label{sec:poisson}
In this section, we briefly discuss the Poisson structure underlying
our spatial discretisation, which will explain the origin of the
conservation of energy and enstrophy. For any functional ${F(\u,h)}$,
${F\colon\mathrm{S}\times\mathrm{V}\to\mathbb{R}}$, we calculate
\begin{align}
\total{F}{t}
  &= \left\langle \dd{F}{\u},\u_t \right\rangle
  + \left\langle \dd{F}{h}, h_t \right\rangle\ ,
\end{align}
where ${\dd{F}{\u}\in S}$ satisfies
\begin{equation}
\left\langle \dd{F}{\u}, \w\right\rangle = \lim_{\varepsilon\to 0}
  \frac{1}{\varepsilon}\left(F(\u+\varepsilon \w,h)-F(\u,h)\right)\ ,
  \quad\forall \w \in \mathrm{S}\ ,
\end{equation}
and similarly ${\dd{F}{h}\in V}$ satisfies
\begin{equation}
\left\langle \dd{F}{h}, \phi\right\rangle = \lim_{\varepsilon\to 0}
  \frac{1}{\varepsilon}\left(F(\u,h+\varepsilon\phi)-F(\u,h)\right)\ ,
  \quad\forall \phi \in \mathrm{V}\ .
\end{equation}
Proceeding with the calculation, we obtain
\begin{align}
\nonumber 
\frac{\mathrm{d} F}{\mathrm{d} t}
  &= \left\langle \dd{F}{\u},-q\F^\perp\right\rangle
  + \left\langle \nabla \cdot \dd{F}{\u}, gh + \frac{|\u|^2}{2} \right\rangle 
  - \left\langle \dd{F}{h}, \nabla\cdot\F \right\rangle \\
&= \left\langle \dd{F}{\u},-q \dd{H}{\u}^\perp\right\rangle
  + \left\langle \nabla \cdot \dd{F}{\u}, \dd{H}{h} \right\rangle
  - \left\langle \dd{F}{h}, \nabla\cdot\dd{H}{\u} \right\rangle
:= \left\{F,H\right\}\ ,
\label{eq:bracket}
\end{align}
where $H$ is the Hamiltonian defined by
\begin{equation}
H = \frac{1}{2}\left\langle \u,h\u \right\rangle
  + \frac{1}{2}\left\langle gh,h \right\rangle\ .
\end{equation}
Equation \eqref{eq:bracket} defines a bilinear bracket for functions
${\mathrm{S}\times\mathrm{V}\to\mathbb{R}}$, which is antisymmetric by
inspection.  This bracket is the restriction to finite elements of a
standard Poisson bracket for shallow-water dynamics. Since we have not
proven the Jacobi identity for the finite element bracket, we only
know that it is an almost-Poisson bracket.

We obtain energy conservation immediately, since ${\dot{H}=\{H,H\}=0}$.
It turns out that enstrophy ${C=\langle q, qh\rangle}$ is a Casimir
for this bracket, since
\begin{align}
\nonumber
\delta C &= \left\langle 2 \delta q, qh \right\rangle
  + \left\langle q^2,\delta h\right\rangle \\
&= \left\langle 2\nabla^\perp\delta q, \u \right\rangle
  + \left\langle 2\delta q, f \right\rangle
  + \left\langle q^2,\delta h \right\rangle\ , 
\end{align}
and therefore ${\dd{C}{\u}=-2\nabla^\perp q}$
(since ${\nabla^\perp q\in \mathrm{S}}$), and
\begin{align}
\left\langle \dd{C}{h}, \delta h \right\rangle
  = \left\langle q^2, \delta h\right\rangle\ .
  \quad\forall \delta h \in V. 
\end{align}
Hence, for any functional $G$,
\begin{align}
\{C,G\} & = \left\langle 2q\nabla^\perp q, \dd{G}{\u}^\perp\right\rangle
  + \left\langle\underbrace{\nabla\cdot -2\nabla^\perp q}_{=0}, \dd{G}{h}\right\rangle
  - \left\langle q^2, \nabla\cdot\dd{G}{\u} \right \rangle \\
&= \left\langle \nabla q^2, \dd{G}{\u}\right \rangle
  - \left\langle q^2,\nabla\cdot\dd{G}{\u} \right \rangle \\
&= 0\ ,\nonumber
\end{align}
where we may integrate by parts in the last line since ${q\in \mathrm{E}}$ and
${\u\in \mathrm{S}}$. $C$ vanishes in the bracket with any other functional and
therefore is a Casimir, i.e.\ a conserved quantity for any
choice of $H$. Unfortunately, there are no known Poisson time
integrators for this type of nonlinear bracket; in particular, the
implicit midpoint rule is not a Poisson integrator for this bracket.

\section*{Acknowledgements}
Andrew McRae wishes to acknowledge funding and other support from the
Grantham Institute. Colin Cotter wishes to acknowledge funding from
NERC grants NE/I02013X/1, NE/I000747/1 and NE/I016007/1.

\bibliography{paper}

\begin{thebibliography}{43}
\providecommand{\natexlab}[1]{#1}
\providecommand{\url}[1]{\texttt{#1}}
\expandafter\ifx\csname urlstyle\endcsname\relax
  \providecommand{\doi}[1]{doi: #1}\else
  \providecommand{\doi}{doi: \begingroup \urlstyle{rm}\Url}\fi

\bibitem[Arakawa(1966)]{arakawa1966computational}
Akio Arakawa.
\newblock Computational design for long-term numerical integration of the
  equations of fluid motion: {T}wo-dimensional incompressible flow. {P}art {I}.
\newblock \emph{Journal of Computational Physics}, 1\penalty0 (1):\penalty0
  119--143, 1966.
\newblock \doi{10.1016/0021-9991(66)90015-5}.

\bibitem[Arakawa and Hsu(1990)]{arakawa1990energy}
Akio Arakawa and Yueh-Jiuan~G. Hsu.
\newblock Energy conserving and potential-enstrophy dissipating schemes for the
  shallow water equations.
\newblock \emph{Monthly Weather Review}, 118\penalty0 (10):\penalty0
  1960--1969, 1990.
\newblock \doi{10.1175/1520-0493(1990)118$<$1960:ECAPED$>$2.0.CO;2}.

\bibitem[Arakawa and Lamb(1977)]{arakawa1977computational}
Akio Arakawa and Vivian~R. Lamb.
\newblock Computational design of the basic dynamical processes of the {UCLA}
  general circulation model.
\newblock In \emph{General Circulation Models of the Atmosphere}, volume~17 of
  \emph{Methods in Computational Physics: Advances in Research and
  Applications}, pages 173--265. Elsevier, 1977.
\newblock \doi{10.1016/B978-0-12-460817-7.50009-4}.

\bibitem[Arakawa and Lamb(1981)]{arakawa1981potential}
Akio Arakawa and Vivian~R. Lamb.
\newblock A potential enstrophy and energy conserving scheme for the shallow
  water equations.
\newblock \emph{Monthly Weather Review}, 109\penalty0 (1):\penalty0 18--36,
  1981.
\newblock \doi{10.1175/1520-0493(1981)109$<$0018:APEAEC$>$2.0.CO;2}.

\bibitem[Arnold et~al.(2006)Arnold, Falk, and Winther]{arnold2006finite}
Douglas~N. Arnold, Richard~S. Falk, and Ragnar Winther.
\newblock Finite element exterior calculus, homological techniques, and
  applications.
\newblock \emph{Acta Numerica}, 15:\penalty0 1--155, 2006.
\newblock \doi{10.1017/S0962492906210018}.

\bibitem[Arnold et~al.(2010)Arnold, Falk, and Winther]{arnold2010finite}
Douglas~N. Arnold, Richard~S. Falk, and Ragnar Winther.
\newblock Finite element exterior calculus: from {H}odge theory to numerical
  stability.
\newblock \emph{Bulletin (New Series) of the American Mathematical Society},
  47\penalty0 (2):\penalty0 281--354, 2010.
\newblock \doi{10.1090/S0273-0979-10-01278-4}.

\bibitem[Bonaventura and Ringler(2005)]{bonaventura2005analysis}
Luca Bonaventura and Todd Ringler.
\newblock Analysis of discrete shallow-water models on geodesic {D}elaunay
  grids with {C}-type staggering.
\newblock \emph{Monthly Weather Review}, 133\penalty0 (8):\penalty0 2351--2373,
  2005.
\newblock \doi{10.1175/MWR2986.1}.

\bibitem[Brezzi and Fortin(1991)]{brezzi1991mixed}
Franco Brezzi and Michel Fortin.
\newblock \emph{Mixed and Hybrid Finite Element Methods}.
\newblock Springer Series in Computational Mathematics. Springer-Verlag, 1991.
\newblock ISBN 0-387-97582-9.

\bibitem[Brezzi et~al.(1985)Brezzi, Douglas, and Marini]{brezzi1985two}
Franco Brezzi, Jim Douglas, Jr., and L.~D. Marini.
\newblock Two families of mixed finite elements for second order elliptic
  problems.
\newblock \emph{Numerische Mathematik}, 47\penalty0 (2):\penalty0 217--235,
  1985.
\newblock \doi{10.1007/BF01389710}.

\bibitem[Brooks and Hughes(1982)]{brooks1982streamline}
AN~Brooks and Thomas~JR Hughes.
\newblock Streamline upwind/{P}etrov--{G}alerkin formulations for convection
  dominated flows with particular emphasis on the incompressible
  {N}avier--{S}tokes equations.
\newblock \emph{Computer Methods in Applied Mechanics and Engineering},
  32\penalty0 (1--3):\penalty0 199--259, 1982.
\newblock \doi{10.1016/0045-7825(82)90071-8}.

\bibitem[Comblen et~al.(2010)Comblen, Lambrechts, Remacle, and
  Legat]{comblen2010practical}
Richard Comblen, Jonathan Lambrechts, Jean-Fran{\c{c}}ois Remacle, and Vincent
  Legat.
\newblock Practical evaluation of five partly discontinuous finite element
  pairs for the non-conservative shallow water equations.
\newblock \emph{International Journal for Numerical Methods in Fluids},
  63\penalty0 (6):\penalty0 701--724, 2010.
\newblock \doi{10.1002/fld.2094}.

\bibitem[Cotter and Ham(2011)]{cotter2011numerical}
C.~J. Cotter and D.~A. Ham.
\newblock Numerical wave propagation for the triangular
  $\mathrm{P}_1^\mathrm{DG}$--$\mathrm{P}_2$ finite element pair.
\newblock \emph{Journal of Computational Physics}, 230\penalty0 (8):\penalty0
  2806--2820, 2011.
\newblock \doi{10.1016/j.jcp.2010.12.024}.

\bibitem[Cotter and Shipton(2012)]{cotter2012mixed}
C.~J. Cotter and J.~Shipton.
\newblock Mixed finite elements for numerical weather prediction.
\newblock \emph{Journal of Computational Physics}, 231:\penalty0 7076--7091,
  2012.
\newblock \doi{10.1016/j.jcp.2012.05.020}.

\bibitem[Danilov(2010)]{danilov2010utility}
Sergey Danilov.
\newblock On utility of triangular {C}-grid type discretization for numerical
  modeling of large-scale ocean flows.
\newblock \emph{Ocean Dynamics}, 60\penalty0 (6):\penalty0 1361--1369, 2010.
\newblock \doi{10.1007/s10236-010-0339-6}.

\bibitem[Danilov et~al.(2008)Danilov, Wang, Losch, Sidorenko, and
  Schr{\"o}ter]{danilov2008modeling}
Sergey Danilov, Qiang Wang, Martin Losch, Dmitry Sidorenko, and Jens
  Schr{\"o}ter.
\newblock Modeling ocean circulation on unstructured meshes: comparison of two
  horizontal discretizations.
\newblock \emph{Ocean Dynamics}, 58\penalty0 (5-6):\penalty0 365--374, 2008.
\newblock \doi{10.1007/s10236-008-0138-5}.

\bibitem[Donea(1984)]{donea1984taylor}
Jean Donea.
\newblock A {T}aylor--{G}alerkin method for convective transport problems.
\newblock \emph{International Journal for Numerical Methods in Engineering},
  20\penalty0 (1):\penalty0 101--119, 1984.
\newblock \doi{10.1002/nme.1620200108}.

\bibitem[Gassmann(2011)]{gassmann2011inspection}
Almut Gassmann.
\newblock Inspection of hexagonal and triangular c-grid discretizations of the
  shallow water equations.
\newblock \emph{Journal of Computational Physics}, 230\penalty0 (7):\penalty0
  2706--2721, 2011.
\newblock \doi{10.1016/j.jcp.2011.01.014}.

\bibitem[Gassmann and Herzog(2008)]{gassmann2008towards}
Almut Gassmann and Hans-Joachim Herzog.
\newblock Towards a consistent numerical compressible non-hydrostatic model
  using generalized {H}amiltonian tools.
\newblock \emph{Quarterly Journal of the Royal Meteorological Society},
  134\penalty0 (635):\penalty0 1597--1613, 2008.
\newblock \doi{10.1002/qj.297}.

\bibitem[Geuzaine and Remacle(2009)]{geuzaine2009gmsh}
Christophe Geuzaine and Jean-Fran{\c{c}}ois Remacle.
\newblock Gmsh: A 3-{D} finite element mesh generator with built-in pre-and
  post-processing facilities.
\newblock \emph{International Journal for Numerical Methods in Engineering},
  79\penalty0 (11):\penalty0 1309--1331, 2009.
\newblock \doi{10.1002/nme.2579}.

\bibitem[Gresho and Sani(1998)]{gresho1998incompressible}
Philip~M Gresho and Robert~L Sani.
\newblock \emph{{Incompressible flow and the finite element method. Volume 1:
  Advection-diffusion and isothermal laminar flow}}.
\newblock John Wiley and Sons, 1998.
\newblock ISBN 0-471-49249-3.

\bibitem[Hallberg and Rhines(1996)]{hallberg1996buoyancy}
Robert Hallberg and Peter Rhines.
\newblock Buoyancy-driven circulation in an ocean basin with isopycnals
  intersecting the sloping boundary.
\newblock \emph{Journal of Physical Oceanography}, 26\penalty0 (6):\penalty0
  913--940, 1996.
\newblock \doi{10.1175/1520-0485(1996)026$<$0913:BDCIAO$>$2.0.CO;2}.

\bibitem[Hyman and Shashkov(1997)]{hyman1997natural}
J.~M. Hyman and M.~Shashkov.
\newblock Natural discretizations for the divergence, gradient, and curl on
  logically rectangular grids.
\newblock \emph{Computers \& Mathematics with Applications}, 33\penalty0
  (4):\penalty0 81--104, 1997.
\newblock \doi{10.1016/S0898-1221(97)00009-6}.

\bibitem[{Le Roux} et~al.(2009){Le Roux}, Hanert, Rostand, and
  Pouliot]{le2009impact}
D.~Y. {Le Roux}, E.~Hanert, V.~Rostand, and B.~Pouliot.
\newblock Impact of mass lumping on gravity and rossby waves in 2{D}
  finite-element shallow-water models.
\newblock \emph{International Journal for Numerical Methods in Fluids},
  59\penalty0 (7):\penalty0 767--790, 2009.
\newblock \doi{10.1002/fld.1837}.

\bibitem[{Le Roux}(2005)]{le2005dispersion}
Daniel~Y. {Le Roux}.
\newblock Dispersion relation analysis of the
  $\mathrm{P}^\mathrm{NC}_1-\mathrm{P}_1$ finite-element pair in shallow-water
  models.
\newblock \emph{SIAM Journal on Scientific Computing}, 27\penalty0
  (2):\penalty0 394--414, 2005.
\newblock \doi{10.1137/030602435}.

\bibitem[{Le Roux}(2012)]{le2012spurious}
Daniel~Y {Le Roux}.
\newblock Spurious inertial oscillations in shallow-water models.
\newblock \emph{Journal of Computational Physics}, 231\penalty0 (6):\penalty0
  7959--7987, 2012.
\newblock \doi{10.1016/j.jcp.2012.04.052}.

\bibitem[{Le Roux} and Pouliot(2008)]{le2008analysis}
Daniel~Y. {Le Roux} and Benoit Pouliot.
\newblock Analysis of numerically induced oscillations in two-dimensional
  finite-element shallow-water models part {II}: Free planetary waves.
\newblock \emph{SIAM Journal on Scientific Computing}, 30\penalty0
  (4):\penalty0 1971--1991, 2008.
\newblock \doi{10.1137/070697872}.

\bibitem[{Le Roux} et~al.(2007){Le Roux}, Rostand, and Pouliot]{le2007analysis}
Daniel~Y. {Le Roux}, Virgile Rostand, and Benoit Pouliot.
\newblock Analysis of numerically induced oscillations in 2{D} finite-element
  shallow-water models part {I}: Inertia-gravity waves.
\newblock \emph{SIAM Journal on Scientific Computing}, 29\penalty0
  (1):\penalty0 331--360, 2007.
\newblock \doi{10.1137/060650106}.

\bibitem[Logg et~al.(2012)Logg, Mardal, and Wells]{logg2012automated}
Anders Logg, Kent-Andre Mardal, and Garth~N. Wells.
\newblock \emph{Automated Solution of Differential Equations by the Finite
  Element Method}.
\newblock Springer, 2012.
\newblock ISBN 978-3-642-23098-1.
\newblock \doi{10.1007/978-3-642-23099-8}.

\bibitem[Raviart and Thomas(1977)]{raviart1977mixed}
P.~A. Raviart and J.~M. Thomas.
\newblock A mixed finite element method for 2-nd order elliptic problems.
\newblock In \emph{Mathematical aspects of finite element methods}, pages
  292--315. Springer, 1977.
\newblock \doi{10.1007/BFb0064470}.

\bibitem[Ringler et~al.(2010)Ringler, Thuburn, Klemp, and
  Skamarock]{ringler2010unified}
T.~D. Ringler, J.~Thuburn, J.~B. Klemp, and W.~C. Skamarock.
\newblock A unified approach to energy conservation and potential vorticity
  dynamics for arbitrarily-structured {C}-grids.
\newblock \emph{Journal of Computational Physics}, 229\penalty0 (9):\penalty0
  3065--3090, 2010.
\newblock \doi{10.1016/j.jcp.2009.12.007}.

\bibitem[Rognes et~al.(2009)Rognes, Kirby, and Logg]{rognes2009efficient}
Marie~E Rognes, Robert~C Kirby, and Anders Logg.
\newblock Efficient assembly of {H}(div) and {H}(curl) conforming finite
  elements.
\newblock \emph{SIAM Journal on Scientific Computing}, 31\penalty0
  (6):\penalty0 4130--4151, 2009.
\newblock \doi{10.1137/08073901X}.

\bibitem[Rostand and Le~Roux(2008)]{rostand2008raviart}
V~Rostand and DY~Le~Roux.
\newblock Raviart--{T}homas and {B}rezzi--{D}ouglas--{M}arini finite-element
  approximations of the shallow-water equations.
\newblock \emph{International Journal for Numerical Methods in Fluids},
  57\penalty0 (8):\penalty0 951--976, 2008.
\newblock \doi{10.1002/fld.1668}.

\bibitem[Sadourny(1975)]{sadourny1975dynamics}
Robert Sadourny.
\newblock The dynamics of finite-difference models of the shallow-water
  equations.
\newblock \emph{Journal of the Atmospheric Sciences}, 32\penalty0 (4):\penalty0
  680--689, 1975.
\newblock \doi{10.1175/1520-0469(1975)032$<$0680:TDOFDM$>$2.0.CO;2}.

\bibitem[Sadourny and Basdevant(1985)]{sadourny1985parameterization}
Robert Sadourny and Claude Basdevant.
\newblock Parameterization of subgrid scale barotropic and baroclinic eddies in
  quasi-geostrophic models: Anticipated potential vorticity method.
\newblock \emph{Journal of the Atmospheric Sciences}, 42\penalty0
  (13):\penalty0 1353--1363, 1985.
\newblock \doi{10.1175/1520-0469(1985)042$<$1353:POSSBA$>$2.0.CO;2}.

\bibitem[Salmon(2005)]{salmon2005general}
Rick Salmon.
\newblock A general method for conserving quantities related to potential
  vorticity in numerical models.
\newblock \emph{Nonlinearity}, 18\penalty0 (5):\penalty0 R1, 2005.

\bibitem[Salmon(2007)]{salmon2007general}
Rick Salmon.
\newblock A general method for conserving energy and potential enstrophy in
  shallow-water models.
\newblock \emph{Journal of the Atmospheric Sciences}, 64\penalty0 (2):\penalty0
  515--531, 2007.
\newblock \doi{10.1175/JAS3837.1}.

\bibitem[Sommer and N{\'e}vir(2009)]{sommer2009conservative}
Matthias Sommer and Peter N{\'e}vir.
\newblock A conservative scheme for the shallow-water system on a staggered
  geodesic grid based on a {N}ambu representation.
\newblock \emph{Quarterly Journal of the Royal Meteorological Society},
  135\penalty0 (639):\penalty0 485--494, 2009.
\newblock \doi{10.1002/qj.368}.

\bibitem[Staniforth and Thuburn(2012)]{staniforth2012horizontal}
Andrew Staniforth and John Thuburn.
\newblock Horizontal grids for global weather and climate prediction models: a
  review.
\newblock \emph{Quarterly Journal of the Royal Meteorological Society},
  138\penalty0 (662):\penalty0 1--26, 2012.
\newblock \doi{10.1002/qj.958}.

\bibitem[Staniforth et~al.(2012)Staniforth, Melvin, and
  Cotter]{staniforth2012analysis}
Andrew Staniforth, Thomas Melvin, and Colin Cotter.
\newblock Analysis of a mixed finite-element pair proposed for an atmospheric
  dynamical core.
\newblock \emph{Quarterly Journal of the Royal Meteorological Society}, 2012.
\newblock \doi{10.1002/qj.2028}.

\bibitem[Thuburn(2008)]{thuburn2008numerical}
J.~Thuburn.
\newblock Numerical wave propagation on the hexagonal {C}-grid.
\newblock \emph{Journal of Computational Physics}, 227\penalty0 (11):\penalty0
  5836--5858, 2008.
\newblock \doi{10.1016/j.jcp.2008.02.010}.

\bibitem[Thuburn and Cotter(2012)]{thuburn2012framework}
J.~Thuburn and C.~J. Cotter.
\newblock A framework for mimetic discretization of the rotating shallow-water
  equations on arbitrary polygonal grids.
\newblock \emph{SIAM Journal on Scientific Computing}, 34\penalty0
  (3):\penalty0 B203--B225, 2012.
\newblock \doi{10.1137/110850293}.

\bibitem[Thuburn et~al.(2009)Thuburn, Ringler, Skamarock, and
  Klemp]{thuburn2009numerical}
J.~Thuburn, T.~D. Ringler, W.~C. Skamarock, and J.~B. Klemp.
\newblock Numerical representation of geostrophic modes on arbitrarily
  structured {C}-grids.
\newblock \emph{Journal of Computational Physics}, 228\penalty0 (22):\penalty0
  8321--8335, 2009.
\newblock \doi{10.1016/j.jcp.2009.08.006}.

\bibitem[Wathen(1987)]{wathen1987realistic}
A.~J. Wathen.
\newblock Realistic eigenvalue bounds for the {G}alerkin mass matrix.
\newblock \emph{IMA Journal of Numerical Analysis}, 7\penalty0 (4):\penalty0
  449--457, 1987.
\newblock \doi{10.1093/imanum/7.4.449}.

\end{thebibliography}

\end{document}